\documentclass[12pt,twoside]{amsart}

\usepackage{amsmath}
\usepackage{amssymb}
\usepackage{amscd}

\renewcommand{\bar}{\overline}

\newcommand{\eps}{\varepsilon}

\newcommand{\CC}{\mathbb{C}}

\newcommand{\HH}{\mathbb{H}}

\newcommand{\PP}{\mathbb{P}}
\newcommand{\QQ}{\mathbb{Q}}
\newcommand{\RR}{\mathbb{R}}

\newcommand{\ZZ}{\mathbb{Z}}

\newcommand{\Qp}{\QQ_p}

\newcommand{\Cp}{\CC_p}

\newcommand{\ints}{{\mathcal O}}

\newcommand{\maxid}{{\mathcal M}}

\newcommand{\calA}{{\mathcal A}}

\newcommand{\calD}{{\mathcal D}}

\newcommand{\calP}{{\mathcal P}}

\newcommand{\calR}{{\mathcal R}}

\newcommand{\PK}{\PP^1(K)}

\newcommand{\PKBerk}{\calP^1(K)}

\newcommand{\PGL}{\hbox{\rm PGL}}

\newcommand{\dist}{\hbox{\rm dist}}

\DeclareMathOperator{\charact}{char}

\newcommand{\Dbar}{\bar{D}}
\newcommand{\DBerk}{\calD}
\newcommand{\DbarBerk}{\bar{\calD}}

\newcommand{\dsps}{\displaystyle}

\theoremstyle{plain}
\newtheorem{thm}{Theorem}[section]

\newtheorem{cor}[thm]{Corollary}
\newtheorem{lemma}[thm]{Lemma}

\newtheorem*{mthm}{Main Theorem}
\newtheorem*{thmblank}{Theorem}

\theoremstyle{definition}
\newtheorem{defin}[thm]{Definition}

\newtheorem*{remarks}{Remarks}
\newtheorem{example}[thm]{Example}

\numberwithin{equation}{section}

\title[Non-archimedean islands]{An Ahlfors Islands Theorem for 
	non-archimedean meromorphic functions}
\author{Robert L. Benedetto}
\date{July 7, 2004}
\thanks{The author gratefully acknowledges the support
  of a Miner D.\ Crary Research Fellowship from Amherst College}
\address{Department of Mathematics and Computer Science \\
        Amherst College \\
        Amherst, MA 01002 \\
        USA}
\email{rlb@cs.amherst.edu}
\urladdr{http://www.cs.amherst.edu/\textasciitilde rlb}

\begin{document}

\newcounter{bean}
\newcounter{sheep}

\begin{abstract}
We present a $p$-adic and non-archimdean version of the
Five Islands Theorem for meromorphic functions from
Ahlfors' theory of covering surfaces.  In the non-archimedean
setting, the theorem requires only four islands, with explicit
constants.  We present examples to show that the constants are
sharp and that other hypotheses of the theorem cannot be removed.
This paper extends an earlier theorem of the author for holomorphic
functions.
\end{abstract}

\maketitle

In the 1930s, Ahlfors proposed his theory of covering surfaces \cite{Ah3}
in complex analysis as an analogue of Nevanlinna theory for domains,
rather than for points.  The Ahlfors theory allows for a description
of the mapping properties of complex meromorphic functions with respect
to open subsets of the image.  One of the key theorems in the subject
is the Five Islands Theorem:
\begin{thmblank}
{\rm (Ahlfors' Complex Five Islands Theorem)}
Let $U_1,\ldots, U_5$ be simply connected domains in the Riemann sphere
with mutually disjoint closures.  Then there
is a constant $h=H(U_1,\ldots,U_5)>0$ with the
following property:
Let $f$ be a complex meromorphic function on the disk $|z|<1$,
and suppose that there is some $r\in (0,1)$ with
\begin{equation}
\label{eq:hcond}
S(f,r)\geq h\cdot L(f,r).
\end{equation}
Then there is a simply connected domain $U$ contained in the disk $|z|<R$
such that
$f$ is one-to-one on $U$ and $f(U)=U_i$ for some $1\leq i\leq 5$.
\end{thmblank}
Here, $S(f,r)$ and $L(f,r)$ (the {\em mean covering number}
and {\em relative boundary length} of $f$, respectively)
are certain real quantities describing the image of $f$ on the open
disk $|z|<r$.
By the work of Dufresnoy \cite{Duf},
condition~\eqref{eq:hcond} may be replaced by a condition
of the form $f^{\#}(0)> \tilde{h}$, where $f^{\#}$ is the
spherical derivative of $f$, and $\tilde{h}$ is, like $h$, a
constant which depends only on the domains $U_1,\ldots, U_5$.
Similar results hold for holomorphic functions, with only
three islands $U_i\subseteq\CC$ required.
Recently, Bergweiler \cite{Berg}
proved the Five Islands Theorem without the theory of covering
surfaces by using a lemma of Zalcman \cite{Zal}, some Nevanlinna theory, and
quasiconformal perturbations.
See \cite{Hay}, Chapters~5--6, for more details on the
theory of covering surfaces.

Initially, the Five Islands Theorem was used mainly in
complex function theory.  Then, in 1968, 
Baker \cite{Bak} applied it in the study of complex dynamics
to prove repelling density for entire functions.
That is, he proved that the Julia
set of a complex entire function must be the closure of the
set of repelling periodic points.
(The usual well known proofs of repelling density
for rational functions do not extend to entire functions.)

In this paper, we will consider non-archimedean fields.
Recall that a non-archimedean field is a field
$K$ equipped with a non-trivial absolute value $|\cdot|$ satisfying
the ultrametric triangle inequality
$|x+y| \leq\max\{ |x|, |y| \}$ for all $x,y\in K$.
Standard examples of such fields include the $p$-adic rationals
$\Qp$ and various function fields.  However, $\Qp$ is not
algebraically closed, and so we set the following notation.
\begin{tabbing}
\hspace{0.5in} \= $K$ \= \hspace{0.5in} \=
  a complete, algebraically closed non-archimedean field with \\
  \> \> \> absolute value $|\cdot|$ \\
\> $\ints_K$ \> \> the ring of integers $\{x\in K : |x|\leq 1\}$ of $K$ \\
\> $k$ \> \> the residue field of $K$
\end{tabbing}
For example, $K$ could be $\Cp$, the completion of an algebraic
closure of $\Qp$.
Recall that the residue field $k$ is defined to
be $\ints_K/\maxid_K$, where
$\maxid_K$ is the maximal ideal $\{x\in K : |x|< 1\}$ of $\ints_K$.
We refer the reader to \cite{Esc,Rob} for treatises on non-archimedean
analysis.

There have been numerous studies
in recent decades of non-archimedean versions of
Nevanlinna theory.  In 1971,
Adams and Straus \cite{AS} proved some
non-archimedean Nevanlinna-style results using
methods much simpler than a full Nevanlinna theory.
More recently, a number of authors have developed
a broader non-archimedean Nevanlinna theory, including analogues
of the First and Second Main Theorems; see \cite{Che1} or
\cite{HuYa} for expositions, and \cite{Bou, BE, CY, Cor, KQ, Ru}
for some of the original papers.

At the same time, there has also been a growing interest in
the dynamics of non-archimedean rational and entire
functions. Broad surveys can be found in
\cite{Ben2, Ben5, Riv1, Riv2}.  Although many
of the fundamental results of complex dynamics have
analogues in the non-archimedean setting, the question
of non-archimedean repelling density remains open,
even in the case of rational functions.
There have been some partial results: Hsia \cite{Hs2}
has shown that the Julia set of a rational function
is contained in the closue of {\em all} periodic points,
and B\'{e}zivin \cite{Bez} has shown that repelling density
follows if there is at least one repelling periodic fixed
point.  However, as discussed in the introduction to \cite{Ben7},
there are serious obstacles to extending either result to
prove repelling density completely.

Bearing Baker's complex result on repelling density in mind,
as well as following the
lead of the non-archimedean Nevanlinna theorists, the author
presented a non-archimedean 
version of Ahlfors' Islands Theorem for holomorphic functions
in \cite{Ben7}.  In that case, only two islands, rather than three,
were required.  However, an extra hypothesis was also needed,
essentially stating that the analogue of $L(f,r)$ is at some point
larger than a constant which depends on the two islands.
In this paper, we continue those investigations by presenting
an analogue of the Islands Theorem for meromorphic functions
in Theorem~\ref{thm:ahlfors}.  We envisage that these
non-archimedean islands theorems should be part of
a non-archimedean theory of Ahlfors' covering surfaces
which is yet to be developed.

The aforementioned
Theorem~\ref{thm:ahlfors} requires the theory
of Berkovich spaces, including the Berkovich projective
line $\PKBerk$,
which we shall discuss in Section~\ref{sect:berk}.  We give
an abbreviated statement of the result here.
\begin{mthm}
{\rm (Non-archimedean Meromorphic Four Islands Theorem)} \\
Let $U_1,U_2,U_3,U_4 \subseteq K\cup\{\infty\}$ be
four disjoint open disks.
Let $\nu_1$ be a Berkovich point such that no connected
component of $\PKBerk \setminus \{ \nu_1\}$ intersects
more than two of $U_1,U_2,U_3,U_4$.
Then there are real constants $C_1,C_2$ depending only on
$K$ and $U_1,U_2,U_3,U_4$ with following property.

Let $f$ be a meromorphic function on $\{z\in K : |z|<1\}$ such that
$f^{\#}(0)>C_1$ and, for any point $\nu\in\DBerk(0,1)$ in the open
Berkovich disk such that $f_*(\nu)=\nu_1$, we have
$\dsps L(f,\nu) \geq C_2$.

Then there is an open disk $U\subseteq D(0,1)$ such that $f$ is one-to-one
on $U$ and $f(U)=U_i$ for some $i=1,2,3,4$.
\end{mthm}
Here, $f^{\#}$ is a non-archimedean version of the spherical
derivative (see Definition~\ref{def:sphderiv} and
equation~\eqref{eq:sphdnu}),
and $L(f,\nu)$ is an analogue of the relative boundary
length (see equation~\eqref{eq:Ldef}).  The other specialized notation,
such as $\DBerk(0,1)$, $\PKBerk$, and $f_*(\nu)$, will be
defined in Section~\ref{sect:berk}.
The full statement of the result, Theorem~\ref{thm:ahlfors},
includes precise descriptions of the constants $C_1$ and $C_2$;
the sharpness of the statement and the constants will be considered
in Examples~\ref{ex:ellip} and~\ref{ex:warp}.

In Section~\ref{sect:proj}, we will recall some basic facts about
the non-archimedean projective line $\PK$.
In Section~\ref{sect:merom} we will review
some standard results about meromorphic functions
on non-archimedean disks.  Section~\ref{sect:berk}
is a summary of the fundamentals of the Berkovich theory.  It is
meant to be a self-contained introduction to the subject, for readers
not familiar with it; proofs of the basic facts will be omitted.
Section~\ref{sect:berk} also includes a number of Lemmas which will
be needed for our main results.
In Section~\ref{sect:bfuncs}, we will study some particular functions
from Berkovich spaces to $\RR$, including the quantity $L(f,\nu)$.
Section~\ref{sect:thm} is devoted to the statement and proof of
the main theorem and a corollary.  Finally, we will present some
examples and address the sharpness of Theorem~\ref{thm:ahlfors}
in Section~\ref{sect:ex}.

The author would like to thank William Cherry
for fruitful discussions on non-archimedean Nevanlinna
and Ahlfors theory, and Jonathan Lubin for some clarification
concerning Example~\ref{ex:ellip}.

\section{The non-archimedean projective line}
\label{sect:proj}

Let $\PK$ denote the projective line over $K$, with
points represented in homogeneous coordinates by $[x,y]$,
for $(x,y)\in K\times K \setminus \{(0,0\}$.  We will
usually identify $\PK$ with $K\cup\{\infty\}$ by taking
$[x,y]$ to $z=x/y$, with $[1,0]$ corresponding to $z=\infty$.

The metric on $K$ induces a standard spherical metric
on $\PK$, given by
$$
\Delta(P_1,P_2) = \frac{|x_1 y_2 - x_2 y_1|}{\max\{|x_1|,|y_1|\}
\max\{|x_2|,|y_2|\} },
$$
where $P_i=[x_i,y_i]$.
Clearly $0\leq\Delta(P_1,P_2)\leq 1$.
In affine coordinates,
$$
\Delta(z_1,z_2) = \frac{|z_1 - z_2|}{\max\{1,|z_1|\} \max\{1,|z_2|\} } .
$$
Note that for $z_1,z_2\in\ints$, we have $\Delta(z_1,z_2)=|z_1-z_2|$.
The topology on $K$ induced by $\Delta$ is exactly the same
as that induced by $|\cdot|$.

The group $\PGL(2,K)$ acts by linear fractional transformations
on $\PK$.  As on the Riemann sphere, given any six points
$P_1,P_2,P_3,Q_1,Q_2,Q_3\in\PK$, there is a unique $\eta\in\PGL(2,K)$
such that $\eta(P_i)=Q_i$ for all $i=1,2,3$.  Of course,
this map $\eta$ need not preserve distances.

On the other hand, the subgroup $\PGL(2,\ints)$ of transformations
$z\mapsto (az+b)/(cz+d)$ with $a,b,c,d\in\ints$ and $|ad-bc|=1$
is distance-preserving with respect to $\Delta$.
(See \cite{Ben7}, Section~1, for example.)  That is,
given $\eta\in\PGL(2,\ints)$ and $P_1,P_2\in\PK$, we have
$$\Delta(\eta(P_1),\eta(P_2)) = \Delta(P_1,P_2).$$
It is easy to check that given any two points $P_1,P_2\in\PK$,
there is a distance-preserving map $\eta\in\PGL(2,\ints)$
such that $\eta(P_1)=P_2$; in fact, there are many such maps.

\section{Holomorphic and meromorphic functions}
\label{sect:merom}

For $a\in K$ and $r>0$, we will denote by $D(a,r)$ and
$\Dbar(a,r)$ the open disk and closed disk (respectively)
of radius $r$ about $a$.  If $r\in|K^{\times}|$,
then $D(a,r)\subsetneq \Dbar(a,r)$,
whereas the two sets coincide if $r\not\in |K^{\times}|$.
It is well known that all disks in $K$
are both open and closed as topological
sets, but we keep the labels ``open disk'' and ``closed disk'' because
the two can behave differently under the action of holomorphic and
meromorphic functions.

By ultrametricity, any point of a disk is a center; but because
$K$ is algebraically closed, the radius is well defined.  That
is, $D(a,r)=D(b,s)$ if and only if $r=s$ and $b\in D(a,r)$;
the analogous statement also holds for closed disks.

\begin{defin}
\label{def:holom}
Let $U\subseteq K$ be a disk.
\begin{list}{\rm \alph{bean}.}{\usecounter{bean}}
\item
  Let $a\in U$ and let $g:U\rightarrow K$.
  We say $g$ is {\em holomorphic} on $U$ if we can write $g$ as a power
  series
  $$g(z) = \sum_{i=0}^{\infty} c_i (z-a)^i \in K[[z-a]]$$
  which converges for all $z\in U$.
\item
  Let $f:U\rightarrow\PK$.
  We say $f$ is {\em meromorphic on $U$} if $f$ is continuous
  on $U$, and if we can write $f$ in homogeneous
  coordinates as
  $$f(z) = [g(z),h(z)]$$
  for all $z$ in some dense subset of $U$,
  where $g$ and $h$ are holomorphic on $U$.
\end{list}
\end{defin}

Thus, a holomorphic function is not just locally analytic but rigid
analytic, in that its defining power series converges on the whole
disk.  (All the rigid analysis in this paper will be hidden from
view inside the results of \cite{Ben7} and \cite{Berk}, but we
refer the interested reader to \cite{BGR,Goss} for detailed background
on the subject.)
As noted in \cite{Ben7}, Section~2, holomorphicity is well defined,
in the sense that if $a,b\in U$, and if $g$ can be written as a convergent
power series centered at $a$, then $g$ can also be written as
a convergent power series centered at $b$.
Naturally, any holomorphic function $f$ is also
meromorphic, by choosing $g=f$ and $h=1$.
Conversely, any meromorphic function which never takes on
the value $\infty$ is in fact holomorphic.

Intuitively, a meromorphic function is simply the
quotient of two holomorphic functions, as in complex analysis.
The technical ``dense subset'' condition in Definition~\ref{def:holom}
is required only because the holomorphic functions $g$ and $h$
may have common zeros.  As observed in \cite{Laz}, it may not
be possible to choose $g$ and $h$ to remove all common zeros if
$U$ is an open disk.  Fortunately, this technicality will not
affect us, because we will not be concerned with any specific
representation $g/h$ of a given meromorphic function $f$.

Derivatives of holomorphic and meromorphic functions are defined
in the same way as in real and complex analysis, and they
satisfy all the usual algebraic rules.  In particular, if $f(z)$
is holomorphic or meromorphic on a disk $U$, then so
is $f'(z)$.

If $f$ is holomorphic on a disk $U$, and if $U'\subsetneq U$
is a smaller disk, then $f(U')$ is a disk.  Moreover,
$f(U')$ is open (respectively, closed) if and only if $U'$ is
open (respectively, closed); see, for example,
\cite{Ben7}, Lemma~2.2.  The following lemma relates the radius
of $f(U')$ to that of $f(U)$.

\begin{lemma}
\label{lem:holom}
Let $a\in K$ and $r>0$.
Let $f$ be a holomorphic function on
the open disk $D(a,r)$
Then
$$
D(f(a), r\cdot |f'(a)| ) \subseteq f(D(a,r))
$$
with equality if $f$ is one-to-one.  The analogous result
also holds for the closed disk $\Dbar(a,r)$.
\end{lemma}

{\bf Proof.}
This follows immediately from \cite{Ben7}, Lemma~2.2.
\qed

As the holomorphic image of a disk is a disk, we are motivated
to define a disk $V\subseteq\PK$ to be either a disk in $K$
in the usual sense or the complement of a disk in $K$.  (Thus,
a disk containing $\infty$ is precisely the complement of
a disk in $K$.)  Equivalently, $V$ is a disk in $\PK$ if
and only if it is the image under some $\eta\in\PGL(2,K)$
of a disk in $K$.  We say $V$ is open (respectively, closed)
if it is either an open (respectively, closed) disk in $K$
or the complement of a closed (respectively, open) disk in $K$.

Given those definitions, if $f$ is meromorphic on a disk $U\subseteq K$,
and if $U'\subsetneq U$ is a strictly smaller disk, then $f(U')$
is either all of $\PK$ or else a disk $V$ in $\PK$.  Moreover,
$V$ is open (respectively, closed) if and only if $U'$ is.

Note that a disk in $\PK$ is not the same as a set of the form
\begin{equation}
\label{eq:badd}
\{ P \in \PK : \Delta(P,a)\leq r\} \quad\text{or}\quad
\{ P \in \PK : \Delta(P,a)<r\}.
\end{equation}
Indeed, for any $r>1$, the region $D(0,r)$ is a disk in $\PK$,
but it cannot be written in the form of \eqref{eq:badd}.  The
same is true of $\PK\setminus D(0,r)$ for any $r<1$.

Even though the spherical metric is not appropriate for defining
radii of disks, it can be made useful for defining derivatives,
as follows.

\begin{defin}
\label{def:sphderiv}
Let $U\subset K$ be a disk, and let $f:U\rightarrow \PK$ be a meromorphic
function.  Let $a\in U$.  The {\em spherical derivative of $f$ at $a$}
is
$$
f^{\#}(a) = \begin{cases}
\dfrac{|f'(a)|}{\max\{1, |f(a)|^2\} } & \text{if } f(a)\neq\infty, \\ \\
\left(\dfrac{1}{f}\right)^{\#}(a) & \text{if } f(a)=\infty.
\end{cases}
$$
\end{defin}

Note that $f^{\#}$ takes values in $[0,\infty)\subseteq\RR$, not in $K$.
The reader may verify that
$$f^{\#}(a) = \lim_{z\rightarrow a} \frac{\Delta(f(z), f(a))}{|z-a|}.$$
Furthermore, if $\eta\in\PGL(2,\ints)$, then because $\eta$ preserves
$\Delta$, it is immediate that
$$(\eta\circ f)^{\#}(a)=f^{\#}(a).$$

We refer the reader to \cite{Esc,Rob} for more on holomorphic
and meromorphic functions.  For an abbreviated survey including
a number of results relevant to this paper,
see \cite{Ben7}, Section~2.

\section{The Berkovich disk and projective line}
\label{sect:berk}

Our discussion of meromorphic functions on $D(0,1)$ will involve
their action on larger spaces defined by Berkovich.
We refer the reader to his papers, especially
\cite{Berk}, for background on general Berkovich spaces and for
proofs of most of their basic properties.
For our purposes, the
reader may find the exposition in \cite{Rum} more useful, as it is
specific to the case of disks and the projective line, which
are all we need here.  The same space for the projective line
was independently discovered later by Rivera-Letelier \cite{Riv1,Riv2};
the set we will call $\PKBerk$ is called $\HH\cup \PK$ in his
notation.

For $a\in K$ and $r\in |K^{\times}|$,
the Berkovich disk $\DbarBerk(a,r)$ associated to the closed disk
$\Dbar(a,r)$ is defined as follows.  Let $\calA(a,r)$ be the
ring of all holomorphic functions on $\Dbar(a,r)$, with Gauss norm
$\nu(a,r)$ given by
\begin{equation}
\label{eq:berkdef}
\|f\|_{\nu(a,r)} = \max\{|c_i| r^i : i\geq 0\},
\end{equation}
where $f(z)=\sum_{i=0}^{\infty}c_i (z-a)^i$.  Intuitively, $\|f\|_{\nu(a,r)}$
is the generic value of $|f(x)|$ on $\Dbar(a,r)$, in the sense that
most $x\in\Dbar(a,r)$ (i.e., all but those in finitely many open subdisks
$D(b,r)$) satisfy $|f(x)|=\|f\|_{\nu(a,r)}$.

A {\em bounded multiplicative seminorm} $\nu$ on $\calA(a,r)$ is a function
$\|\cdot\|_{\nu}:\calA(a,r)\rightarrow [0,\infty)$ such that
for all $f,g\in\calA(a,r)$,
\begin{list}{\rm \roman{bean}.}{\usecounter{bean}}
\item $\|0\|_{\nu} = 0$,
\item $\|1\|_{\nu} = 1$,
\item $\|fg\|_{\nu} = \|f\|_{\nu} \cdot \|g\|_{\nu}$,
\item $\|f+g\|_{\nu} \leq\max\{ \|f\|_{\nu} , \|g\|_{\nu} \}$, and
\item $\|f\|_{\nu} \leq \|f\|_{\nu(a,r)}.$
\end{list}
(The above versions of properties (iv) and (v) are stronger than
the usual definitions, but they are equivalent for our ring
$\calA(a,r)$, as shown in the first few pages of \cite{Rum}.)
The function $\nu$ is called a seminorm because $\|f\|_{\nu}=0$
need not necessarily imply that $f=0$.

The Berkovich disk $\DbarBerk(a,r)$
is then defined to be the set of all bounded multiplicative
seminorms on $\calA(a,r)$, with the Gel'fond topology, which
is the weakest topology such that
$$\{\nu\in \DbarBerk(a,r) : \|f\|_{\nu} < R\}
\quad\text{and}\quad
\{\nu\in \DbarBerk(a,r) : \|f\|_{\nu} > R\}$$
are open, for any $f\in\calA(a,r)$ and any $R\in\RR$.
The reader may check that for fixed $f,g\in\calA(a,r)$
and $R\in\RR$, the set
$$\{\nu\in \DbarBerk(a,r) : \|f\|_{\nu} < R\|g\|_\nu\}$$
is also open.
The space $\DbarBerk(a,r)$ is compact, Hausdorff, and path-connected
\cite{Berk}.  Berkovich also showed that the points of $\DbarBerk(a,r)$
come in four types, which we now list.

There is a natural inclusion of $\Dbar(a,r)$ in $\DbarBerk(a,r)$,
as follows.  If $x\in\Dbar(a,r)\subseteq K$, then the value of
the seminorm $\|f\|_{x}$ is defined to be simply $|f(x)|$.
The corresponding points
of $\DbarBerk(a,r)$ are called the {\em type~I points}.

Meanwhile, for every $b\in\Dbar(a,r)$ and every $s\in(0,r]$,
the seminorm $\nu(b,s)$, which in this case is actually a norm,
is defined exactly as in
equation~\eqref{eq:berkdef}, but for the disk $\Dbar(b,s)$.
That is,
$$
\|f\|_{\nu(b,s)} = \max\{ |d_i| s^i : i\geq 0\},
\quad
\text{where}
\quad
f(z) = \sum_{i=0}^{\infty} d_i (z-b)^i.
$$
Note that this definition makes sense for any $s\in (0,r]$,
not just $s\in |K^{\times}|$.
Also note that the norm depends only on the disk $\Dbar(b,s)$,
not on the choice $b$ of center.
The corresponding point $\nu(b,s)$ of $\DbarBerk(a,r)$
is said to be of {\em type~II} if $s\in |K^{\times}|$ and of
{\em type~III} otherwise.

From the type~II and~III points, one can begin to see how
$\DbarBerk(a,r)$ is path-connected, at least between type~I
points.  Given $b,c\in\Dbar(a,r)$, the path from $b$ to $c$
starts at $b$, which we consider as a disk of radius zero.
We increase the radius through a path of type~II and~III
points of the form $\nu(b,s)$ until we get to $s=|b-c|$.
Then $\nu(b,s)=\nu(c,s)$, and so we may decrease the radius
$s$ towards the new center $c$ until we arrive at $c$ itself.

Finally, there is one more class of points.
A {\em type~IV} point $\nu$
corresponds to a nested infinite sequence of
disks $U_1\supset U_2\supset \cdots$ with empty intersection.
(Such sequences may exist because $K$ need not be maximally
complete.  The infimum of the radii of the $U_i$ for such
a sequence is always strictly positive.)  The norm $\nu$
is simply the limit of the norms $\nu(U_i)$.  (Of course,
it is the limit norm $\nu$, not the sequence $\{\nu(U_i)\}$,
which is the type IV point; there are infinitely many
equivalent sequences
$\{U_i\}$ that approach any given type IV point.)
The type~IV
points are needed to make $\DbarBerk(a,r)$ compact, but they
will not be important in our discussions.

Intuitively,
the space $\DbarBerk(a,r)$ looks like a tree branching out
from the root point $\nu(a,r)$ with infinitely many
branches at every type~II point (which are dense in the tree),
and with limbs ending at the type~I and type~IV points.
The infinitely many branches at a type~II point $\nu(b,s)$
correspond to the infinitely many open subdisks $D(c,s)$ of
$\Dbar(b,s)$ of the same radius, as well as (if $s<r$) one
more branch corresponding to increasing the radius (i.e.,
corresponding to the disk at $\infty$).
The type~III points, meanwhile, are interior points with
no branching.

If $f$ is meromorphic on $\Dbar(a,r)$ and $\nu\in\DbarBerk(a,r)$,
then we may define $\|f\|_{\nu}$ to be $\|g\|_{\nu}/\|h\|_{\nu}$,
where $f=g/h$ for $g,h\in\calA(a,r)$.  (Note that
$\|f\|_{\nu}=\infty$ if and only if $\nu=b$ is type~I and
$f$ has a pole at $b$.)  As before, it is appropriate to
think of $\|f\|_{\nu(a,r)}$ as the generic value of $|f(x)|$
for $x\in\Dbar(a,r)$.
The extended function $\|\cdot\|_{\nu}$
still satisfies properties (i)--(iv) of multiplicative
seminorms, but we may no longer have $\|f\|_{\nu}\leq\|f\|_{\nu(a,r)}$.

For an open disk $D(a,r)$, we can also associate a Berkovich
space $\DBerk(a,r)$ by taking the union (really,
the direct limit) of sets $\DbarBerk(a,r_i)$, where $r_i\nearrow r$.
The resulting space is still path-connected, Hausdorff, and
locally compact, but it is no longer compact.
Although $\DBerk(0,1)$ will be one of our main objects of study,
we will understand it by considering the subspaces $\DbarBerk(a,r)$
described above, for $a\in D(0,1)$ and $0<r<1$.

We may also define the Berkovich projective line $\PKBerk$
by glueing two copies of $\DBerk(0,r)$ (for some $r>1$)
as follows.
A type~I point $x$ on one copy with
$1/r<|x|<r$ is identified with $1/x$ on the other copy.
Meanwhile, a type~II
or~III point $\nu(b,s)$ with $1/r<|b|<r$ is identified with $\nu(1/b,s/|b|^2)$,
since $\Dbar(1/b,s/|b|^2)$ is the image of $\Dbar(b,s)$ under
$z\mapsto 1/z$.  A type~IV point which is the limit of a sequence
of type~II points is mapped to the limit of the image of the sequence
under $z\mapsto 1/z$.

Thus, $\PKBerk$ looks like $\DbarBerk(0,1)$ with an extra copy of
the open tree $\DBerk(0,1)$ attached to the top (i.e., the $\infty$ end)
of the point $\nu(0,1)$.  The new top portion contains all points $x$
of $\PK$ with $|x|>1$, including $\infty$, as well as points
$\nu(a,r)$ with $|a|>1$ or $r>1$.
Like $\DbarBerk(0,1)$, the space $\PKBerk$ is path-connected,
Hausdorff, and compact.

Any disk in $\PK$ is associated with a unique point
(of type~II or~III) of $\PKBerk$.  Indeed, any open or closed
disk $\Dbar(a,r)$ or $D(a,r)$ or its complement
is associated with the point $\nu(a,r)$.
Conversely, a type~III point $\nu(a,r)$ is associated with exactly
two disks, namely $D(a,r)=\Dbar(a,r)$ and its complement.
Meanwhile, a type~II point $\nu(a,r)$ is associated with infinitely
many disks: every open disk $D(b,r)$ for $b\in\Dbar(a,r)$, the
disk $\PK\setminus\Dbar(a,r)$, and the complements of all these
open disks.  Although a type~II or~III point is associated with
more than one disk, note that it is associated with exactly one
closed disk which does not contain $\infty$.  Thus, there is
a one-to-one correspondence between type~II and~III points of $\PKBerk$
and closed disks in $K$; the point $\nu(a,r)$ corresponds to
$\Dbar(a,r)$.

Borrowing from \cite{Riv2}, we state the following definition.
\begin{defin}
\label{def:sep}
Let $X$ be a connected Berkovich space, let $W\subseteq X$ be a subset,
and let $\nu_1\in X$ be a point.  We say that {\em $\nu_1$ separates $W$}
if $W$ intersects more than one connected component of
$X\setminus\{\nu_1\}$.
\end{defin}
Following the intuition of the tree structure,
it is not difficult to show for $X=\PKBerk$, $X=\DBerk(a,r)$,
or $X=\DbarBerk(a,r)$,
that if $\nu_1\in X$
separates any set, then $\nu_1$ must be type~II or~III.

We will usually consider separation in the case that
the subset $W$ contains only type~I points.
For example, $\nu(0,1)$ separates any subset of $\PK$
that contains both $0$ and a point $a$ with $|a|=1$.  However,
$\nu(0,1)$ does not separate $D(0,1)$.  Clearly, if $W\subseteq V$
and $\nu_1$ separates $W$, then $\nu_1$ also separates $V$.

As  mentioned at the start of this section, a meromorphic function
$f$ on $\Dbar(a,r)$ induces a function
$f_*:\DbarBerk(a,r)\rightarrow\PKBerk$.
A fully rigorous derivation of $f_*$ and its properties requires
a description of general Berkovich spaces as locally ringed spaces
with patches given by general Berkovich affinoids.  We refer
the reader to \cite{Berk} or \cite{Rum}, Section~B, for such
a derivation; an equally rigorous derivation in a different style
appears in \cite{Riv2}, Section~4.
We will now describe $f_*$ precisely, but we will skip the proofs.

Using the more general Berkovich machinery, one can show that
for each $\nu\in\PK$, there is a corresponding local ring $\calA_{\nu}$
of functions $f$ for which $\|f\|_{\nu}$ can be defined.  If $\nu$
is of type~II, III, or~IV, then $\calA_{\nu}$ contains $K(z)$,
the ring of rational functions over $K$.  If $\nu=b$ is of type~I,
then $\calA_{\nu}$ contains all functions in $K(z)$ except those
with poles at $b$; in that case, we may still talk about $\|f\|_{b}$
for functions $f$ with a pole at $b$ by defining $\|f\|_{b}=\infty$.
The crucial fact from the Berkovich machinery is that
$\nu$ is completely determined by its restriction to $K(z)$.

Thus, given $f$ meromorphic on $\Dbar(a,r)$ and $\nu\in\DBerk(a,r)$,
we define $f_*(\nu)$ to be the unique point
(i.e., seminorm) in $\PKBerk$ such
that for all $h\in K(z)$,
\begin{equation}
\label{eq:berkmapdef}
\|h\|_{f_*(\nu)} = \| f\circ h \|_{\nu}.
\end{equation}
The same definition applies to a meromorphic function $f$
on $D(a,r)$.
Similarly, if $\eta\in\PGL(2,K)$ and $\nu\in\PKBerk$,
we define $\eta_*(\nu)$
so that $\|h\|_{\eta_*(\nu)} = \| \eta\circ h \|_{\nu}$
for all $h\in K(z)$.

Equation~\eqref{eq:berkmapdef} is somewhat unsatisfying at
first; besides the fact that we have omitted the proof of
existence and uniqueness of $f_*(\nu)$, the equation does
not give much immediate insight into what $f_*$ really
looks like.  Following \cite{Riv2}, then, we present
the following equivalent description.

Let $\Dbar(a,r)$ be a closed disk with $r\in |K^{\times}|$.
If $f$ is holomorphic on $\Dbar(a,r)$
and $\nu\in\DbarBerk(a,r)$ is a point of type~II or~III,
then write $\nu=\nu(b,s)$ for some disk
$D(b,s)\subsetneq \Dbar(a,r)$.  From Section~\ref{sect:merom}
we know that $f(D(b,s))$ is an open disk; write
$f(D(b,s))=D(f(b),\sigma)$.
Then $f_*$ from equation~\eqref{eq:berkmapdef} satisfies
$$f_*(\nu(b,s))=\nu(f(b),\sigma).$$
More generally,
given $f$ meromorphic on $\Dbar(a,r)$ and $\nu\in\DbarBerk(a,r)$
of type~II or~III, write $\nu=\nu(b,s)$ for some disk
$D(b,s)\subsetneq D(a,r)$.
It can be shown that there is a radius $s_0<s$ such that
for all $s'$ with $s_0<s'<s$, the image
$f(D(b,s)\setminus \Dbar(b,s'))$ of the
annulus $D(b,s)\setminus \Dbar(b,s')$ is itself an annulus
of the form
$$D(\beta,\sigma)\setminus \Dbar(\beta,\sigma')
\quad\text{or}\quad
D(\beta,\sigma')\setminus \Dbar(\beta,\sigma),$$
where $\beta$ and $\sigma$ are fixed and do not vary with
$s'$.  Then it turns out that
$$f_*(\nu(b,s))=\nu(\beta,\sigma).$$

It follows quickly from the definitions that $f_*$
is a continuous function from one Berkovich space to another,
and that $f_*$ agrees with $f$ at the type~I points.
Moreover, if $f$ is a nonconstant meromorphic function,
then $f_*$ takes type~I points to type~I points,
type~II points to type~II points, and so on.

Given $f$ meromorphic on $\Dbar(a,r)$, we can extend
$f^{\#}$ from the type~I points to all of $\DbarBerk(a,r)$
by setting
\begin{equation}
\label{eq:sphdnu}
f^{\#}(\nu) = \frac{\|f'\|_{\nu}}{\max\{ 1, \|f\|_{\nu}^2 \} }.
\end{equation}
As was true for Definition~\ref{def:sphderiv}, it is easy to
show that if $\eta\in\PGL(2,\ints)$, then
$(\eta\circ f)^{\#}(\nu)=f^{\#}(\nu)$.

\begin{defin}
\label{def:sphrad}
Given a point $\nu(a,\rho)\in \PKBerk$ of type~II or~III, define
\begin{equation}
\label{eq:rdef}
r(\nu(a,\rho)) = \frac{\rho}{\max\{1, \rho^2, |a|^2\} }
\end{equation}
to be the {\em spherical radius} of $\nu(a,\rho)$.
\end{defin}

We leave it to the reader to verify that $r$ is independent of
the choice of $a$ in $\Dbar(a,\rho)$, and that $r$
is a continuous function 
on the subset of $\PKBerk$ on which it is defined.
The spherical radius may also be defined at points of type~I and~IV
by continuity, so that $r:\PKBerk\rightarrow [0,\infty)$ is continuous.
In fact, $r(\nu)$ is just $\exp[-\dist(\nu, \nu(0,1))]$,
where $\dist(\cdot,\cdot)$ is the metric on $\PKBerk\setminus\PK$
which appears in \cite{Riv2}, Section~3 as the metric on $\HH$
and in \cite{Rum}, Section~B as the ``big model'' metric.

If we restrict $r(\cdot)$ to $\DBerk(0,1)$ or $\DbarBerk(0,1)$,
then $r(\nu(a,\rho))=\rho$ is the usual radius of the associated
disk.  Similarly, if $X=\PKBerk$ and if $\nu=\nu(a,\rho)$ separates
$\Dbar(0,1)$ (which is to say that $\Dbar(a,\rho)\subseteq\Dbar(0,1)$),
then $r(\nu)$ is again just the usual radius $\rho$.  On the other
hand, if $\nu(a,\rho)$ does not separate $\Dbar(0,1)$, then $r(\nu)<\rho$.
Formula \eqref{eq:rdef} is chosen so that $r$
is invariant under distance-preserving transformations, as
the following lemma shows.

\begin{lemma}
\label{lem:etar}
Let $\nu\in\PKBerk$, and let $\eta\in\PGL(2,\ints)$.
Then $r(\eta_*(\nu)) = r(\nu)$.
\end{lemma}

{\bf Sketch of Proof.}
We may assume $\nu$ is not of type I or IV, as the result
for types II and III will extend by continuity.
Since $\PGL(2,\ints)$ is generated by maps of the
form $z+A$ (for $|A|\leq 1$), $Bz$ (for $|B|=1$), and $1/z$,
we may consider only such maps.  The verification
is trivial for $z+A$ and $Bz$.  For $\eta(z)=1/z$, we
may write $\nu=\nu(a,\rho)$ with $\rho>0$.
If $0\in\Dbar(a,\rho)$, then
$\nu=\nu(0,\rho)$ and $\eta_*(\nu) = \nu(0,1/\rho)$,
from which the verification of the lemma is easy.
On the other hand, if $0\not\in\Dbar(a,\rho)$, then
$|a|>\rho$, and $\eta_*(\nu) = \nu(1/a, \rho/|a|^2)$.
The lemma then follows.
\qed

The remaining more specific lemmas will be needed to prove
Theorem~\ref{thm:ahlfors}.

\begin{lemma}
\label{lem:berkext}
Let $f$ be a nonconstant meromorphic function on $D(0,1)$,
let $x\in D(0,1)$, let $0<r<1$,
and set $\nu_1 = f_*(\nu(x,r))\in\PKBerk$.
Then $\nu_1$ separates $f(D(x,r+\eps))$ for every
$\eps>0$.
\end{lemma}

{\bf Proof.}
Replacing $f$ by $\eta\circ f$ for an appropriate
$\eta\in\PGL(2,K)$, we may assume that $f(x)=0$
and $\nu_1=\nu(0,\rho)$ for some $\rho>0$.  If
$\nu_1$ separates $f(D(x,r))$, then we are done.
If not, then $f(x)=0$ forces $f(D(x,r))\subseteq D(0,\rho)$;
because $f_*(\nu(x,r))=\nu_1$, we must have
$f(D(x,r))=D(0,\rho)$.  For any given $\eps>0$,
if $f$ has a pole in $D(x,r+\eps)$, then we are
done; so we may assume that $f$ is holomorphic on
$D(x,r+\eps)$.  Since $f$ is nonconstant with $f(x)=0$,
$\|f\|_{\nu(x,s)}$ must be a strictly increasing function
of $s$ for $0<s<r+\eps$; see equation~\eqref{eq:berkdef}.
It follows that $\|f\|_{\nu(x,r+\eps)}>\rho$, and therefore
$f(D(x,r+\eps))$ contains a point $a$ with $|a|>\rho$,
implying that $\nu_1$ separates $f(D(x,r+\eps))$.
\qed

\begin{lemma}
\label{lem:berkfin}
Let $f$ be meromorphic on $D(0,1)$,
let $x\in D(0,1)$, let $0<r<1$,
and let $\nu_1\in\PKBerk$.  Then
there are only finitely many points
$\nu\in\DBerk(x,r)$ such that
$f_*(\nu)=\nu_1$.
\end{lemma}

{\bf Proof.}
This Lemma follows easily from the machinery of
\cite{Riv2}, Section~4, but we include a direct
proof for the convenience of the reader.

Write $f=g/h$ for $g$ and $h$ holomorphic.
If $\nu_1=a\in\PK$ is a type I point, then by a change
of coordinates, we may assume $a=0$.
The lemma then holds for $\nu_1$ by
the finiteness of the Weierstrass degree of $g$
on $D(x,r)\subsetneq D(0,1)$; see, for example,
\cite{Ben7}, Lemma~2.2.

If $\nu_1$ is of type II or III, then let $a=f(x)$; there must
be a point $b\in\PK$ such that $\nu_1$ separates $\{a,b\}$.
By a change of coordinates, we may assume $a=0$ and $b=\infty$,
so that $\nu_1=\nu(0,R)$ for some $R>0$.
By the previous paragraph, there are only finitely
many zeros and poles of $f$ in $D(x,r)$.

If $f_*(\nu(y,s))=\nu_1$
for some $\Dbar(y,s)\subseteq D(x,r)$, then
because $f(\Dbar(y,s))$ is either $\PK$ or a closed disk
associated with $\nu_1$, there must be a root of
$f=0$ or $f=\infty$ in $\Dbar(y,s)$.
Thus, there can be only finitely many disjoint disks
$\Dbar(y_i,s_i)$ with $f_*(\nu(y_i,s_i))=\nu_1$, or
else there would be infinitely many poles or zeros of $f$
in $D(x,r)$, contradicting the previous paragraph.

Meanwhile, if $y\in D(x,r)$ and $0<s_1<s_2< r$
with $f_*(\nu(y,s_i))=\nu_1$ for $i=1,2$, we claim
that there must be a root $y'$ of $f=0$ or $f=\infty$
with $s_1 < |y'-y| \leq s_2$.  If not, then
move $y$ to $0$ and
write the holomorphic functions $g$, $h$ as
$g(z) = \sum_{j=0}^{\infty} a_j z^j$
and
$h(z) = \sum_{j=0}^{\infty} b_j z^j$.
Our assumption about the lack of roots implies that
one term $a_m z^m$ of $g$ and one term $b_n z^n$
of $h$ is uniquely
maximal in each sum for all $s_1 < |z| \leq s_2$.
Thus, for all such $z$,
$$|f(z) - c z^{\ell}| < |f(z)|,$$
where $c=a_m/b_n\in K$ and $\ell=m-n\in\ZZ$.
Since $f_*(\nu(y,s_1))=f_*(\nu(y,s_2))=\nu(0,R)$,
we must have $|c z^{\ell}|=R$ for all such $z$,
which implies that $\ell=0$ and $|c|=R$.
In that case, however, $|f(z)-c|<R$ for all such $z$,
meaning in particular that $|f(z)-c|<R$ for
all $|z|=s_2$.  Then $f_*(\nu(0,s_2))\neq \nu(0,R)$,
which is a contradiction and proves our claim.

Thus, any chain $\Dbar(y_1,s_1)\supsetneq \Dbar(y_2,s_2)\supsetneq\cdots$
of disks with $f_*(\nu(y_i,s_i))=\nu_1$ must be finite,
or else there would be infinitely many poles or zeros of $f$
in $D(x,r)$.
Together with the above fact that only finitely many
disjoint disks $\Dbar(y_i,s_i)$ can have $f_*(\nu(y_i,s_i))=\nu_1$,
it follows that there can be only
finitely many $\nu\in\DBerk(x,r)$ such that $f_*(\nu)=\nu_1$,
as desired.

We will not need to consider the case that $\nu_1$
is type IV in this paper, and we leave the proof
of that case to the reader.
\qed

\begin{lemma}
\label{lem:berkptmap}
Let $f$ be meromorphic on $D(x,r)$ for some $x\in K$
and $r>0$, and let $\nu_1\in\PKBerk$.  Then
there is a small enough radius $r'>0$ such that
$\nu_1$ does not separate $f(D(x,r'))$.
\end{lemma}

{\bf Proof.}
By a change of coordinates, we may assume that $f(x)=0$.
The statement is vacuous if $\nu_1$ is type I or IV,
so we assume it is type II or III.  We may therefore
write $\nu_1=\nu(a,R)$ for some $a\in K$, $R>0$.

By Lemma~\ref{lem:berkfin}, $f$ has only finitely many poles in
$D(x,r/2)$.  We may therefore choose $0<s <r/2$ such
that there are no poles in $D(x,s)$, implying that
$f$ is holomorphic on $D(x,s)$.
If $f(D(x,s))\subseteq D(0,R)$, we are done.
Otherwise, by \cite{Ben7}, Lemma~2.6, there is a radius $r'\in (0,s]$
such that $f(D(x,r'))=D(0,R)$, which is
not separated by $\nu_1$.
\qed

\begin{lemma}
\label{lem:berkinf}
Let $f$ be meromorphic on $D(x,r)$ for some $x\in K$
and $r>0$, and let $\nu_1\in\PKBerk$.  Suppose that
$\nu_1$ separates $f(D(x,r))$.  Then
\begin{list}{\rm \alph{bean}.}{\usecounter{bean}}
\item there is some $\eps>0$ such that $\nu_1$ separates
  $f(D(x,r-\eps))$, and
\item there is a closed disk $\Dbar(z,s)\subseteq D(x,r)$
  such that $f_*(\nu(z,s))=\nu_1$.
\end{list}
\end{lemma}

{\bf Proof.}
Since $\nu_1$ separates $f(D(x,r))$, there must be a point
$y\in D(x,r)$ such that $\nu_1$ separates $\{f(x),f(y)\}$.
Let $\eps=(r-|x-y|)/2>0$.  Then $x,y\in D(x,r-\eps)$,
so that $\nu_1$ separates $f(D(x,r-\eps))$, proving part (a).

Recall that $\DBerk(x,r)$ is connected and that
$f_*:\DBerk(x,r) \rightarrow \PKBerk$ is continuous.
Because $x,y\in\DBerk(x,r)$ but $f(x)$ and $f(y)$
are in different connected components of $\PKBerk\setminus\{\nu_1\}$,
there must be some $\nu\in\DBerk(x,r)$ such that
$f_*(\nu)=\nu_1$.  Note that $\nu_1$, and therefore $\nu$,
must be points of type II or III, because $\nu_1$ separates
$\PK$.  Thus, we may write $\nu=\nu(z,s)$ for some
closed disk $\Dbar(z,s)\subseteq D(x,r)$.
\qed

\section{Functions on Berkovich spaces}
\label{sect:bfuncs}

Theorem~\ref{thm:ahlfors} and its proof will rely heavily
on certain functions from $\DBerk(0,1)$ to $\RR$.  We have
already seen the spherical radius function $r(\nu)$ defined
in equation~\eqref{eq:rdef}.  In addition,
if $f$ is a holomorphic function on $D(0,1)$, then $f$ induces
another real-valued map, given by
$\nu \mapsto \|f\|_{\nu}$.
Such maps are continuous, because if $f=g/h$ with $g,h$ holomorphic,
then
$$\left\{\nu : \left\| \frac{g}{h} \right\|_{\nu} <R \right\}
= \{ \nu: \|g\|_{\nu} < R \|h\|_{\nu} \}$$
is open in the Gel'fond topology.  The following lemma
gives a more precise description of the behavior of $\|f\|_{\nu}$.

\begin{lemma}
\label{lem:monom}
Let $f$ be a nonconstant meromorphic function on $D(0,1)$, and let
$a\in D(0,1)$.  For $0<r<1$, define $F(r) = \|f\|_{\nu(a,r)}$.
Then
\begin{list}{\rm \alph{bean}.}{\usecounter{bean}}
\item $F:(0,1)\rightarrow (0,\infty)$ is a continuous function
  which is piecewise of the form $F(r)=c r^n$, for $c\in (0,\infty)$
  and $n\in\ZZ$.
\item If $F(r)=c r^n$ on the interval
  $[r_1,r_2)$ for some $0<r_1<r_2\leq 1$, then
  $n$ is the number of zeros of $f$ in $\Dbar(a,r_1)$ less the
  number of poles, counting multiplicity of each.
\item For any fixed $0<R<1$, there are only finitely many radii
  $r\in (0,R]$ at which the value of the exponent $n$ can change.
\end{list}
\end{lemma}

{\bf Proof.}
This is a combination of Lemmas~4.4 and~4.5 of \cite{Ben7}.
\qed

We now use the map $\nu\mapsto\|f\|_{\nu}$ to construct
several other more specialized functions for use in our main theorem.
Given a meromorphic function $f$ on $D(0,1)$,
define
$L:\DBerk(0,1) \rightarrow [0,\infty)$ by
\begin{equation}
\label{eq:Ldef}
L(\nu) = L(f,\nu) = r(\nu) \cdot f^{\#}(\nu)
= \frac{r(\nu) \cdot \|f'\|_{\nu}}{\max\{1,\|f\|_{\nu}\}}.
\end{equation}
We use the notation $L$ because the above function is
a non-archimedean analogue of Ahlfors' relative boundary length
function $L(f,r)$.
Note that for $\eta\in\PGL(2,\ints)$, we have
$L(f,\nu) = L(\eta\circ f, \nu)$, because $f^{\#} = (\eta\circ f)^{\#}$.

Next, given $f$ and a point $\alpha\in\Dbar(0,1)$
with $\alpha\neq 0,1$, define
$G:\DBerk(0,1) \rightarrow [0,\infty)$ by
\begin{equation}
\label{eq:Gdef}
G(\nu) = G(f,\alpha,\nu) = \frac{\left(r(\nu) \cdot \|f'\|_{\nu}\right)^2}
{\|f\|_{\nu} \|f-\alpha\|_{\nu} \|f- 1\|_{\nu}}.
\end{equation}
The reason for the conditions on $\alpha$ will become clear
in the proof of Theorem~\ref{thm:ahlfors}.

Equation~\eqref{eq:Gdef} currently does not make sense
if $\nu$ is one of the four points $0,\alpha,1,\infty$ of type~I.
However, $G$ extends naturally
to all of $\DBerk(0,1)$, as we now argue.
By Lemma~\ref{lem:monom}, for any fixed $a\in D(0,1)$,
$G(\nu(a,r))$ is a continuous, piecewise
monomial function of $0<r<1$, as is $L(\nu(a,r))$.
In fact, by the same Lemma, for
any $\|f\|_{\nu(a,r)}$, there is some $\eps>0$ such that
the exponent $n$ in $F(r)=cr^n$ is constant on $0<r<\eps$ and equal
to the order of the zero (or negative the order of the pole)
of $f$ at $a$.  It then follows fairly easily that the definition
of $G$ extends to all points of $\DBerk(0,1)$,
so that $G$ is continuous and finite-valued on $\DBerk(0,1)$.

Finally, given $f$ and a point $b\in\PK$, then
for any closed disk $\Dbar(a,r)\subseteq D(0,1)$,
define
\begin{equation}
\label{eq:Ndef}
\begin{split}
  N_b(a,r) &= \text{number of roots of $f=b$ in $\Dbar(a,r)$,} \,
     \; \text{and} \\
  N_{\text{ram}}(a,r) &= \text{number of ramification points of $f$
    in $\Dbar(a,r)$,}
\end{split}
\end{equation}
where in each case ``number'' means the number of points, counted
with multiplicity.
Note that $N_{\text{ram}}$ counts with multiplicity all points
at which $f'=0$; but it also counts ramification at
all multiple poles.  More precisely,
if $f=g/h$ with $g$ and $h$ holomorphic, then $N_{\text{ram}}$ counts
(with multipicity) the zeros of $g'h-h'g$, less twice the number
of common zeros of $g$ and $h$.
In addition, note that by Lemma~\ref{lem:monom}, $G(\nu(a,r))$ is locally
a monomial in $r$ of degree
$$2 + 2 N_{\text{ram}} - (N_0 + N_{\alpha} + N_1 + N_{\infty}),$$
because we may write
$$G(\nu) = \frac{\left(r(\nu) \cdot \|g'h-h'g\|_{\nu}\right)^2}
{\|g\|_{\nu} \|g-\alpha h \|_{\nu} \|g- h\|_{\nu} \|h\|_{\nu} } .$$

We now list several properties of $L$, $G$, and $N$
which will be useful in the proof of Theorem~\ref{thm:ahlfors}.

\begin{lemma}
\label{lem:Linv}
Let $f$ be a meromorphic function on $D(0,1)$, let
$\eta\in\PGL(2,K)$, and let $\nu\in\DBerk(0,1)$ be a point
of type II or III.  Then
$$\frac{L(f,\nu)}{r(f_*(\nu))} =
\frac{L(\eta\circ f , \nu)}{r(\eta_*(f_*(\nu_1)))} .$$
\end{lemma}

{\bf Proof.}
Let $\nu_1=f_*(\nu)$, and
write $\nu_1=\nu(a,\rho)$, with $a\in K$ and $\rho>0$.
Then $\|f\|_{\nu} = \max\{|a|,\rho\}$, so that
\begin{equation}
\label{eq:Ltemp}
\frac{L(f,\nu)}{r(f_*(\nu))}
=
\frac{r(\nu) \|f'\|_{\nu} }{r(\nu_1) \max\{1, \|f\|_{\nu}^2\} }
= 
\frac{r(\nu) \|f'\|_{\nu} \max\{1,\rho^2, |a|^2\} }
{\rho \max\{1, \rho^2, |a|^2\}}
=
\frac{r(\nu) \|f'\|_{\nu}}{\rho} .
\end{equation}

We may factor $\eta=\eta_1\circ \eta_2$, where
$\eta_1\in\PGL(2,\ints)$ and $\eta_2(\infty)=\infty$.
Since $\eta_1$ preserves $r$ and $L$,
we may assume without loss that $\eta=\eta_2$.
Thus, $\eta(z)=B(z-a) + A$ for some $A, B\in K$ with $B\neq 0$.

We compute $\|(\eta\circ f)'\|_{\nu} = |B| \cdot \|f'\|_{\nu}$
(by the chain rule) and $\eta_*(\nu_1)=\nu(A, |B|\rho)$.
Thus, by equation~\eqref{eq:Ltemp},
$$
\frac{L(\eta\circ f,\nu)}{r(\eta_*(f_*(\nu)))}
=
\frac{r(\nu) \|(\eta\circ f)'\|_{\nu}}{|B| \rho}
=
\frac{r(\nu) \|f'\|_{\nu}}{\rho}
=
\frac{L(f,\nu)}{r(f_*(\nu))}.
$$
\qed

The quantity $L(f,\nu)/r(f_*(\nu))$ in Lemma~\ref{lem:Linv}
is a measure of distortion; it generalizes the quantity
$\delta$ which appeared in \cite{Ben7}.  Intuitively, $L$ itself measures
the expected spherical radius of $f_*(\nu)$ based on the generic
value of $f^{\#}$ at $\nu$.  Because $f^{\#}$ may be smaller
than one would suspect (since $z^p$ has small derivative
$pz^{p-1}$ in residue characteristic $p$), the actual spherical
radius $r(f_*(\nu))$ of the image may be larger than $L$.
Thus, the smaller the distortion ratio, the further the actual
spherical radius is from that predicted by the derivative.

In this paper, we will only need to consider the
quantity $L(f,\nu)/r(f_*(\nu))$ at points of type II or III.
Nonetheless, we note here (and leave to the reader to verify)
that the same quantity extends by continuity to all of $\DBerk(0,1)$.
If $\nu=x$ is of type~I,
then $L(f,\nu)/r(\nu_1)$ (which is $0/0$ as written)
turns out to be $|n|$, where $n\geq 1$ is the multiplicity
with which $x$ maps to $f(x)$.

The functions $L$ and $G$ are bounded above by $1$,
as the following lemma states.

\begin{lemma}
\label{lem:Gbound}
Let $f$ be a meromorphic function on $D(0,1)$,
and let $\alpha\in\Dbar(0,1)\setminus D(1,1)$, with $\alpha\neq 0$.
Then for any $\nu\in\DBerk(0,1)$,
$$\left( L(f,\nu) \right)^2 \leq G(f,\alpha,\nu) \leq 1.$$
\end{lemma}

{\bf Proof.}
By continuity, it suffices to show the result in the case that $\nu$
is type~II.
Recall from \cite{Ben7}, Lemma~4.2, that
$r(\nu) \|g'\|_{\nu} \leq \|g\|_{\nu}$
for a holomorphic function $g$ on $\DBerk(0,1)$.
Thus, writing $f=g/h$ for $g,h$
holomorphic on $D(0,1)$, we see that
\begin{equation}
\label{eq:delta}
r(\nu) \|f'\|_{\nu} = \frac{r(\nu)\|g'h - h'g\|_{\nu} }{\|h\|_{\nu}^2}
\leq \frac{\|g\|_{\nu}}{\|h\|_{\nu}} \cdot \max \left\{
\frac{r(\nu)\|g'\|_{\nu}}{\|g\|_{\nu} } , 
\frac{r(\nu)\|h'\|_{\nu}}{\|h\|_{\nu} } \right\}
\leq \|f\|_{\nu},
\end{equation}
so that the same inequality holds for meromorphic functions
as well.

Write $f_*(\nu)=\nu(a,R)$.
If $f_*(\nu)$ separates $\PK\setminus\Dbar(0,1)$,
which is to say that $\max\{|a|,R\} > 1$, then
$$\|f\|_{\nu} = \|f-\alpha\|_{\nu} = \|f-1\|_{\nu} > 1,$$
so that by inequality~\eqref{eq:delta},
$$G(\nu)= \frac{r(\nu)^2 \|f'\|_{\nu}^2}{ \|f\|_{\nu}^3}
\leq \frac{1}{\|f\|_{\nu} } < 1.$$
The inequality $L^2\leq G$ follows similarly.

If $f_*(\nu)$ separates $D(0,|\alpha|)$,
which is to say that $\max\{|a|,R\} < |\alpha|$, then
$$\|f\|_{\nu} < \|f-\alpha\|_{\nu} = |\alpha|,
\quad\text{and}\quad
\|f-1\|_{\nu} = 1,$$
so that by inequality~\eqref{eq:delta},
$$G(\nu)=\frac{r(\nu)^2 \|f'\|_{\nu}^2}{|\alpha| \|f\|_{\nu}}
< \left(\frac{r(\nu) \|f'\|_{\nu} }{\|f\|_{\nu} } \right)^2 \leq 1;$$
again, the $L^2\leq G$ inequality also follows easily.

If $f_*(\nu)$ separates either $D(\alpha,|\alpha|)$
or $D(1,1)$, the verification is similar.  Thus, the
only remaining case is that
$\alpha$ and $R$ satisfy
$|\alpha|\leq \max\{|a|,R\}\leq 1$, 
$|\alpha|\leq \max\{|a-\alpha|,R\}\leq 1$, 
and
$\max\{|a-1|,R\}= 1$.
In that case,
$$\|f\|_{\nu} =\|f-\alpha\|_{\nu} \leq 1
\quad\text{and}\quad
\|f-1\|_{\nu}=1,$$
so that
$$G(\nu)=\frac{r(\nu)^2 \|f'\|_{\nu}^2}{\|f\|_{\nu}^2} \leq 1,$$
and $L(\nu)^2=r(\nu)^2 \|f'\|_{\nu}^2 \leq G(\nu)$.
\qed

The following lemma, which uses the notation
of equation~\eqref{eq:Ndef}, will be useful
for choosing more useful centers for certain disks.

\begin{lemma}
\label{lem:center}
Let $f$ be a meromorphic function on $D(0,1)$,
and let $\alpha\in\Dbar(0,1)\setminus\{0,1\}$.
Fix $x\in D(0,1)$ and a radius $0<R<1$, and suppose that
$$N_0(x,R) + N_{\alpha}(x,R) + N_1(x,R) + N_{\infty}(x,R)
  > 2 N_{\text{ram}}(x,R).$$
Then there is a point $y\in \Dbar(x,R)$ such that for every
$r\in [0,R]$,
$$N_0(y,r) + N_{\alpha}(y,r) + N_1(y,r) + N_{\infty}(y,r)
  > 2 N_{\text{ram}}(y,r).$$
\end{lemma}

{\bf Proof.}
Suppose not.  Define
$$N_{\text{tot}}(a,r) = 
N_0(a,r) + N_{\alpha}(a,r) + N_1(a,r) + N_{\infty}(a,r).$$
By Lemma~\ref{lem:berkfin}, there are only
finitely many roots $\{y_i\}_{i=1}^m$ of
$f=0,\alpha,1,\infty$ in $\Dbar(x,R)$.
Then for each $i=1,\ldots, m$, there is some $r_i\in [0,R]$
such that $N_{\text{tot}}(y_i,r_i) \leq 2 N_{\text{ram}}(y_i,r_i)$.
If any two of these disks intersect,
then one contains the other, and so we may discard the smaller
one.  We are left with a finite set $\{\Dbar(y'_i,r'_i)\}_{i=1}^{\ell}$
of pairwise disjoint disks in $\Dbar(x,R)$, each of which
satisfies $N_{\text{tot}}(y'_i,r'_i) \leq 2 N_{\text{ram}}(y'_i,r'_i)$,
and which together contain all of the $\{y_i\}$.  Thus,
$$ N_{\text{tot}}(x,R) = \sum_{i=1}^{\ell} N_{\text{tot}}(y'_i,r'_i)
\leq 2\sum_{i=1}^{\ell} N_{\text{ram}}(y'_i,r'_i)
\leq 2 N_{\text{ram}}(x,R),$$
contradicting the hypotheses and hence proving the lemma.
\qed

\section{The Four Islands Theorem}
\label{sect:thm}

We need the following specialized radius to define
the key value $\mu$ which will appear in Theorem~\ref{thm:ahlfors}.

\begin{defin}
\label{def:ahlrad}
Let $U_1,U_2,U_3,U_4\subseteq\PK$ be four disjoint open disks
in $\PK$, and for each $i=1,2,3,4$, fix a point $a_i\in U_i$.
Choose a map $\eta_1\in\PGL(2,K)$ such that $\eta_1(a_1)=0$,
$\eta_1(a_2)=\infty$, and $\eta_1(a_3)=1$.  Write
$\eta_1(U_1)=D(0,r_1)$, and define
$$s_1 =\frac{r_1}{\min\{|\eta(a_2)|,|\eta(a_3)|,|\eta(a_4)|\} }
=\frac{r_1}{\min\{1,|\eta(a_4)|\} }.$$
Define $s_i$ similarly for $i=2,3,4$ by mapping $a_i$ to $0$
and two of the other points to $\infty$ and $1$.
We define the {\em Ahlfors radius} of $\{U_1,U_2,U_3,U_4\}$
to be
$$s = \max\{s_1,s_2,s_3,s_4\}.$$
\end{defin}

The reader may verify that the Ahlfors radius
(and even the set $\{s_1,s_2,s_3,s_4\}$) is well defined,
in the sense that it is independent of the ordering of the
indices $1,2,3,4$ and the choice of the $\{a_i\}$.
In addition, if $\tilde{\eta}\in\PGL(2,K)$,
then the Ahlfors radius of
$\{\tilde{\eta}(U_1),\tilde{\eta}(U_2),
\tilde{\eta}(U_3),\tilde{\eta}(U_4)\}$ is the same
as that of $\{U_1,U_2,U_3,U_4\}$.
Note that $0<s\leq 1$.

Intuitively,
each $s_i$ is the ratio of the radius of $U_i$ to the distance
from $U_i$ to the nearest other $U_j$.  The Ahlfors radius
is then the largest of these ratios.  In the course of our
proof, we will move the points $\{a_i\}$ to $\{0,\alpha,1,\infty\}$
for some $\alpha\in\Dbar(0,1)\setminus D(1,1)$ with $\alpha\neq 0$.
To say that
the Ahlfors radius of disks centered at those four points is
at most $s$ is to say that the four disks are contained in
$D(0,|\alpha|s)$,
$D(\alpha,|\alpha|s)$, $D(1,s)$, and $\PK\setminus\Dbar(0,1/s)$.

We are now prepared to state our main result.

\begin{thm}
\label{thm:ahlfors}
{\rm (Non-archimedean Meromorphic Four Islands Theorem)} \\
Let $K$ be a complete, algebraically closed non-archimedean
field with residue field $k$, and let $p=\charact k \geq 0$.
If $p\geq 3$, define the real number
$$E_p = \sum_{i=1}^{\infty} \frac{1}{p^i - 1}.$$
Let $U_1,U_2,U_3,U_4 \subseteq \PK$ be
four pairwise disjoint open disks.
Let $\nu_1\in\PKBerk$ such that no connected
component of $\PKBerk \setminus \{ \nu_1\}$ intersects
more than two of $U_1,U_2,U_3,U_4$.

Let $s$ be the Ahlfors radius of $\{U_1, U_2, U_3, U_4\}$.  Set
$$\mu=\left\{
\begin{array}{ll}
0 & \mbox{if $\charact k = 0$,} \\
s^{1/2} & \mbox{if $\charact k = 2$,} \\
\min\left\{s^{(\frac{1}{2}-\frac{1}{2p})}, |p|^{-E_p} s^{1/2}\right\} &
\mbox{if $\charact k = p\geq 3$,}
\end{array}
\right. $$
and set
$$
C_1 = \frac{1}{r(\nu_1)},
\quad \mbox{and}\quad
C_2 = \mu \cdot r(\nu_1).$$

Let $f$ be a meromorphic function in $D(0,1)$ such that
\begin{list}{\rm \alph{bean}.}{\usecounter{bean}}
\item $\dsps f^{\#}(0)>C_1$, and
\item for any point $\nu\in\DBerk(0,1)$
  such that $f_*(\nu)=\nu_1$, we have
  $\dsps L(f,\nu) \geq C_2$.
\end{list}
Then there is an open disk $U\subseteq D(0,1)$ such that $f$ is one-to-one
on $U$ and $f(U)=U_i$ for some $i=1,2,3,4$.
\end{thm}

\begin{remarks}
\hspace{1pt}

\begin{list}{\rm \arabic{bean}.}{\usecounter{bean}}
\item
Because the Ahlfors radius satisfies $0<s\leq 1$, it
follows that $0\leq\mu\leq 1$.

\item
The statement of the theorem becomes stronger if $\mu$
is smaller, because more functions $f$ will satisfy condition~(b).
Informally, then, the smaller $\mu$ is, the better.

\item
Suppose $\charact k=p\geq 3$.  If $\charact K=0$, then $\mu$
decreases on the order of $s^{1/2}$ as $s$ approaches $0$.
On the other hand, if $\charact K =p$, then $\mu$ is
$s^{(1/2 - 1/(2p))}$, which decreases more slowly.

\item
In the $\charact k = 0$ case, the choice of $\mu=0$ above
means that the value of $C_2$, and therefore the value of the
Ahlfors radius $s$, is irrelevant; condition~(b) becomes vacuous.

\item
The lower bound of $C_1$ for $f^{\#}(0)$ in condition~(a)
is sharp, as we now observe.
Choose $\lambda\in K$ with $|\lambda|\geq 1$.
Let $f(z)=\lambda z$, 
let $\nu_1=\nu(0,|\lambda|)$, and
let each $U_i$ be
a disk of the form $D(a_i,\eps)$, with $|a_i|=|\lambda|$ and
$\eps>0$ as small as one wishes.
Note that $f^{\#}(0)=C_1=|\lambda|$, so that condition~(a)
just barely fails; however, all the other hypotheses
of Theorem~\ref{thm:ahlfors} hold.
Nonetheless, the image $f(D(0,1))$
fails to intersect any $U_i$, let alone map a subdomain
onto one of them.

The sharpness of condition~(b) is more subtle and
will be considered in the examples of Section~\ref{sect:ex}.
\end{list}
\end{remarks}

The following Lemma will appear in the final step of the proof
of the Theorem.  It is a slightly more complicated version
of \cite{Ben7}, Proposition~5.2.

\begin{lemma}
\label{lem:mudisk}
Let $K$, $k$, $p$, and $E_p$ be as in Theorem~\ref{thm:ahlfors}.
Let $\alpha\in\Dbar(0,1)\setminus D(1,1)$ with $\alpha\neq 0$,
and let $\nu_1\in\PKBerk$.  Suppose that no connected component of
$\PKBerk \setminus \{ \nu_1\}$ contains more than
two of $0,\alpha,1,\infty$.
Let $f$ be a meromorphic function on $D(0,1)$,
let $\mu\in [0,1]$, let $x\in D(0,1)$, and let $0<R'<1$.
Suppose that
\begin{list}{\rm (\roman{bean})}{\usecounter{bean}}
\item $\nu_1$ separates $f(D(x,R'+\eps))$ for every $\eps>0$,
\item $\nu_1$ does not separate $f(D(x,R'))$,
\item $N_0(x,r) +N_{\alpha}(x,r) + N_1(x,r)
+ N_{\infty}(x,r) > 2 N_{\text{ram}}(x,r)$
for all $r\in [0,R']$, and
\item $G(f,\alpha,\nu(x,R'))\geq \mu^2$.
\end{list}
Then there is an open disk $U\subseteq\Dbar(x,R')$ such that $f$ is
one-to-one on $U$, and $f(U)$ is one of
$$
D(0,|\alpha|s), \quad
D(\alpha,|\alpha|s), \quad
D(1,s), \quad \text{or} \quad
\PK \setminus \Dbar(0,1/s),
$$
where
$$s= \left\{
\begin{array}{ll}
1 & \mbox{if } p = 0 \\
\mu^2 & \mbox{if } p = 2 \\
\max\left\{ |p|^{2E_p} \mu^2, \mu^{2p/(p-1)} \right\} &
\mbox{if } p \geq 3.
\end{array} \right.
$$
\end{lemma}

Condition~(i) is not actually required to prove the Lemma.
It is stated here only for convenience, as all four conditions
(i)--(iv) will figure prominently in the proof of
Theorem~\ref{thm:ahlfors}.
If $\nu_1$ is a type~II Berkovich point, then condition~(i)
is equivalent to the simpler statement that
$\nu_1$ separates $f(\Dbar(0,r_0))$.  However, the more complicated
statement is required if $\nu_1$ is type~III.

{\bf Proof of Lemma~\ref{lem:mudisk}.}
We may write $\nu_1=\nu(0,\rho)$ for some $|\alpha|\leq\rho\leq 1$.
By choosing $r$ sufficiently small in property~(iii),
we have $f(x)\in\{0,\alpha,1,\infty\}$.

We may assume that
$f(x)=0$.  Indeed, if $f(x)=\alpha$, then let $\tilde{f}=\eta\circ f$,
where $\eta(z) = (z-\alpha)/(1-\alpha)$, which takes
$\alpha$ to $0$, fixes $1$ and $\infty$, and takes $0$
to $\tilde{\alpha}=\alpha/(\alpha-1)$.  The corresponding
$\tilde{\nu}_1$ and $G(\tilde{f},\tilde{\alpha},\cdot)$ satisfy 
conditions (i)--(iv), and the disks $\eta(D(0,|\alpha|s))$,
$\eta(D(\alpha,|\alpha|s)$, etc., are $D(\tilde{\alpha}, |\tilde{\alpha}|)$,
$D(0, |\tilde{\alpha}|)$, etc., as appropriate.
Similar arguments hold for $f(x)=\infty$
(with $\eta(z)=\alpha/z$) and for $f(x)=1$ (with
$\eta(z)=[\alpha(z-1)]/[(\alpha-1) z]$).

Since $\nu_1$ does not separate $f(D(x,R'))$, and since $f(x)=0$,
$f$ cannot have poles in $D(x,R')$.  Similarly, it cannot take
on the value $1$ in $D(x,R')$.
Thus,
$f$ and $f'$ are holomorphic on $D(x,R')$, and $N_1(x,r)=N_{\infty}(x,r)=0$
for all $0<r<R'$.

Let $R''=\inf\{0<r\leq R' : \alpha\in f(\Dbar(x,r))\}$ if this
set is nonempty, or $R''=R'$ if it is empty.  For $R''<r<R'$,
we have $N_0(x,r)=N_{\alpha}(x,r)$, because any two points
in the image of a holomorphic function have the same number
of preimages; see \cite{Ben7}, Lemma~2.2.  For any such $r$,
then, condition~(iii) becomes
\begin{equation}
\label{eq:Neven}
2 N_{\text{ram}}(x,r) + 2 -
(N_0(x,r) +N_{\alpha}(x,r) + N_1(x,r) + N_{\infty}(x,r)) \leq 0,
\end{equation}
because $N_0 + N_{\alpha} + N_1 + N_{\infty}=2N_0$ is even.
However, the left side of \eqref{eq:Neven}
is the local monomial degree of $G(\nu(x,r))$ at $r$.
Thus, $G(\nu(x,r))$ is a locally constant or decreasing function
of $r$ on $(R'',R')$.  By continuity and condition~(iv), we have
$$G(\nu(x,R'')) \geq G(\nu(x,R')) \geq \mu^2.$$

Because $f(D(x,R''))\subseteq D(0,|\alpha|)$, we have
$\|f-\alpha\|_{\nu(x,r)} = |\alpha|$ and $\|f-1\|_{\nu(x,r)}= 1$
for all $0<r\leq R''$.  For all such $r$, then,
\begin{equation}
\label{eq:Gmod}
G(\nu(x,r)) =
\frac{r^2 \|f'\|_{\nu(x,r)}^2}{|\alpha| \cdot \|f\|_{\nu(x,r)} } ,
\end{equation}
and for $r<R''$, condition (iii) is
$$N_0(x,r) \geq 1 +  2 N_{\text{ram}}(x,r).$$

The remainder of the proof is identical with that of
\cite{Ben7}, Proposition~5.2,
from the statement of Lemma~5.5
onward (pages~613--615).  If we replace $f$ by $f/\alpha$,
then expression~\eqref{eq:Gmod} becomes 
$r^2 \|f'\|_{\nu(x,r)}^2/\|f\|_{\nu(x,r)}$, which
is exactly $F(r)$ in the notation of \cite{Ben7}.
The proof in \cite{Ben7} ultimately concludes that this new $f$
maps some disk $U=D(x,\tilde{R})$ one-to-one onto $D(0,s)$.
Thus, the original map $f$ maps $U$
one-to-one onto $D(0,|\alpha|s)$, and we are done.
\qed

{\bf Proof of Theorem~\ref{thm:ahlfors}.}
For each $i=1,2,3,4$, choose a point $a_i\in U_i$.
By hypothesis, we may assume that neither $a_1$ nor $a_2$ lies in
the same component of $\PKBerk\setminus\{\nu_1\}$ as either $a_3$ or $a_4$.
Similarly,
by exchanging the roles of $\{a_1,a_2\}$ and $\{a_3,a_4\}$ if necessary,
we may assume that $f(0)$ does not lie in the same component
as either $a_3$ or $a_4$.

Replacing $f$ by $\eta_0\circ f$ for some $\eta_0\in\PGL(2,\ints)$,
we may assume that $a_4=\infty$; this change does not affect
any of the hypotheses, by Section~\ref{sect:proj}, by
Lemma~\ref{lem:Linv}, and by the
discussion following Definition~\ref{def:ahlrad}.
Note that $\nu_1$ separates $a_1$ from both $a_3$ and $a_4=\infty$,
so that $\nu_1=\nu(a_1,\rho)$ for some $0<\rho\leq |a_3-a_1|$.
Also note that
$a_2,f(0)\in \Dbar(a_1,\rho)\setminus D(a_3,|a_3-a_1|)$.
By the hypothesis that $f^{\#}(0)> 1/r(\nu_1)\geq \rho$, we compute
$$|f'(0)| = f^{\#}(0) \cdot \max \{1,|f(0)|^2\} > \rho.$$

Define $\eta\in\PGL(2,K)$ by $\eta(z) = (z-a_1)/(a_3-a_1)$, which is
chosen so that $\eta(a_1)=0$, $\eta(a_3)=1$, and $\eta(a_4)=\infty$.
Write $\tilde{f}=\eta\circ f$,
$\tilde{\rho}=\rho/|a_3-a_1| \leq 1$,
$\tilde{\nu}_1 = \eta_*(\nu_1)$,
and $\alpha=\eta(a_2)$.
Then $\tilde{\nu}_1 = \nu(0,\tilde{\rho})$,
$\alpha\in\Dbar(0,\tilde{\rho})\setminus D(1,1)$,
and $\tilde{f}(0)\in\Dbar(0,\tilde{\rho})$.
In addition,
$$|\tilde{f}'(0)| = \frac{|f'(0)|}{|a_3-a_1|}
> \frac{\rho}{|a_3-a_1|} = \tilde{\rho}.$$

By Lemma~\ref{lem:Linv}, the hypothesis that
$L(f,\nu) \geq \mu r(\nu_1)$ for any $\nu$ with $f_*(\nu)=\nu_1$
becomes $L(\tilde{f},\nu) \geq \mu r(\tilde{\nu}_1) = \mu\tilde{\rho}$
for any $\nu$ with $\tilde{f}_*(\nu)=\tilde{\nu}_1$.
Any such $\nu$ satisfies
$$\| \tilde{f}\|_{\nu} = \| \tilde{f}-\alpha \|_{\nu} = \tilde{\rho},
\quad\text{and}\quad
\| \tilde{f}-1 \|_{\nu} = 1,$$
so that $L(\tilde{f},\nu)=r(\nu) \|\tilde{f}'\|_{\nu}$,
and therefore
$$
G(\tilde{f},\nu)= \frac{r(\nu)^2 \|\tilde{f}'\|_{\nu}^2}
{\| \tilde{f}\|_{\nu} \| \tilde{f}-\alpha \|_{\nu} \| \tilde{f}-1 \|_{\nu} }
= \frac{L(\tilde{f},\nu)^2}{\tilde{\rho}^2}
\geq \mu^2.
$$

Let
$$\calR_0=\left\{ 0<r<1 : \nu_1
\text{ does not separate } \tilde{f}(D(0,r)) \right\},$$
which is nonempty, by Lemma~\ref{lem:berkptmap}.
Let $r_0=\sup\calR_0>0$.  If $r_0=1$, then
$\tilde{f}(D(0,1))\subseteq \Dbar(0,\tilde{\rho})$.
In that case, $\tilde{f}$ is holomorphic on $D(0,1)$,
and by Lemma~\ref{lem:holom}, $|\tilde{f}'(0)|\leq \tilde{\rho}$,
which is a contradiction.  Thus, $0<r_0<1$.

If $\nu_1$ separates $\tilde{f}(D(0,r_0))$, then by
Lemma~\ref{lem:berkinf}.a, $\nu_1$ also separates
$\tilde{f}(D(0,r_0-\eps))$ for some $\eps>0$.  In that
case, $\calR_0\cap [r_0-\eps, r_0]=\emptyset$, which
contradicts the fact that $r_0$ is the supremum of $\calR_0$.
Thus, $\nu_1$ does not separate $\tilde{f}(D(0,r_0))$, but
for every $\eps>0$, $\nu_1$ separates $\tilde{f}(D(0,r_0+\eps))$.
Moreover,
$\tilde{f}(D(0,r_0))\subseteq D(f(0),\tilde{\rho})
\subseteq \Dbar(0,\tilde{\rho})\setminus D(1,1)$,
and $\tilde{f}$ is holomorphic on $(D(0,r_0))$.


From now on, we will no longer need the hypothesis
that $f^{\#}(0)> C_1$.  Writing $f$ in place of $\tilde{f}$, then,
we have $0<r_0<1$, $0 < \rho \leq 1$, $\nu_1=\nu(0,\rho)$,
$\alpha\in\Dbar(0,\rho)\setminus D(1,1)$ with $\alpha\neq 0$,
and a meromorphic function $f$ on $D(0,1)$ such that
\begin{itemize}
\item $f$ is holomorphic on $D(0,r_0)$,
\item $f(D(0,r_0))\subseteq\Dbar(0,\rho)\setminus D(1,1)$,
\item $|f'(0)|>\rho$,
\item $\nu_1$ does not separate $f(D(0,r_0))$,
\item $\nu_1$ separates $f(D(0,r_0+\eps))$ for all $\eps>0$, and
\item for every $\nu\in\DBerk(0,1)$ such that $f_*(\nu)=\nu_1$,
  we have $G(\nu)\geq \mu^2$.
\end{itemize}
We wish to show that
$f$ maps some open disk $U\subseteq D(0,1)$ one-to-one
onto one of $D(0,s/|\alpha|)$, $D(\alpha,s/|\alpha|)$,
$D(1,s)$, or $\PK\setminus\Dbar(0,1/s)$.

Since $f(D(0,r_0))\subseteq\Dbar(0,\rho)\setminus D(1,1)$, we have
$$\| f \|_{\nu(0,r_0)} ,  \| f-\alpha \|_{\nu(0,r_0)} \leq \rho,
\quad\text{and}\quad
\| f-1 \|_{\nu} = 1.$$
Moreover, since $f$ (and hence $f'$) is holomorphic on
$D(0,r_0)$, we have $\|f'\|_{\nu(0,r_0)}\geq |f'(0)|$;
see equation~\eqref{eq:berkdef}.
Thus,
\begin{equation}
\label{eq:Gbig}
G(\nu(0,r_0))= \frac{r_0^2 \|f'\|_{\nu(0,r_0)}^2}
{\| f\|_{\nu(0,r_0)} \| f-\alpha \|_{\nu(0,r_0)} \| f-1 \|_{\nu(0,r_0)} }
\geq \frac{|f'(0)|^2}{\rho^2} \cdot r_0^2.
\end{equation}

Let
$$\calR = \{r\in [r_0,1) :
  N_0(0,r) + N_{\alpha}(0,r) + N_1(0,r) + N_{\infty}(0,r)
  > 2 N_{\text{ram}}(0,r)\}.
$$
As noted in Section~\ref{sect:bfuncs},
$G(\nu(0,r))$ is locally a monomial function of $r$ of degree
$$2 + 2 N_{\text{ram}}(0,r) -
[N_0(0,r) + N_{\alpha}(0,r) + N_1(0,r) + N_{\infty}(0,r) ].$$
Therefore, a radius $r\in [r_0,1)$ is in $\calR$
if and only if $G(\nu(0,\cdot))$ is of degree
strictly less than two at $r$.

We claim that $\calR\neq\emptyset$.  Indeed, if $\calR=\emptyset$,
then by Lemma~\ref{lem:monom}, $G(\nu(0,\cdot))$ is a continuous function
which is piecewise monomial, and always of degree at least two,
on $[r_0,1)$.
By inequality~\eqref{eq:Gbig}, it follows that
$$G(\nu(0,r))\geq\left(\frac{|f'(0)|}{\rho}\right)^2 r^2$$
for all $r\in [r_0,1)$.  Since $|f'(0)|>\rho$, there must be
some such $r$ for which $G(\nu(0,r))>1$, contradicting
Lemma~\ref{lem:Gbound} and proving the claim.

Let $R\in\calR$.  Note that $R\geq r_0$; thus,
$\nu_1$ separates $f(\Dbar(0,R+\eps))$ for every $\eps>0$.
Moreover, by Lemma~\ref{lem:center}, we may choose
$y\in\Dbar(0,R)$ such that for every $r\in [0,R]$,
$$N_0(y,r) + N_{\alpha}(y,r) + N_1(y,r) + N_{\infty}(y,r)
  > 2 N_{\text{ram}}(y,r).$$
We wish to find a (possibly different) point $x\in\Dbar(0,R)$
and a (possibly smaller) radius $R'$ so that the pair $(x,R')$
satisfies properties (i)--(iv) of Lemma~\ref{lem:mudisk}.

We will do so by an inductive process.  We begin
with $y_1=y$ and an auxiliary radius $R_1=R$.
We have just observed that properties (i) and (iii) already
apply to the pair $(y_1,R_1)$.
At each step $n\geq 1$, given a
point $y_n$ and an auxiliary radius $R_n$,
we will define the radius $R'_n$.
We will prove that
all of the desired properties hold for $(y_n,R'_n)$ except
possibly condition~(iv).
If that condition fails, we will construct
a new point $y_{n+1}$ and a new auxiliary radius $R_{n+1}$,
and the process will repeat.
We will then prove that there must eventually
be some $n\geq 1$ for which 
condition~(iv) holds.

The process is as follows.
At step $n\geq 1$, we are given $y_n$ and an auxiliary radius $R_n$
such that properties (i) and (iii) apply to the pair
$(y_n,R_n)$.  Define
$$\calR'_n = \{ r\in (0,R_n] :
\nu_1 \text{ separates } f(D(y_n,r+\eps)) \text{ for every } \eps>0 \},$$
which is nonempty because $R_n\in \calR'_n$.
Let $R'_n=\inf\calR'_n$.
Observe that $R'_n>0$, by Lemma~\ref{lem:berkptmap}.
By definition of $R'_n$ and the properties of $(y_n,R_n)$,
properties (i) and (iii) apply to the pair $(y_n,R'_n)$.
In addition, property (ii) applies to the pair, by
Lemma~\ref{lem:berkinf}.a.
If $G(\nu(y_n,R'_n))\geq \mu^2$, then property~(iv)
holds, and our process finishes by setting $(x,R')=(y_n,R'_n)$.
We may therefore assume that $G(\nu(y_n,R'_n))< \mu^2$.

We claim there is a disk
$\Dbar(z_{n+1},R''_{n+1})\subseteq \Dbar(y_n,R'_n)$
such that $f_*(\nu(z_{n+1},R''_{n+1})) = \nu_1$.
By definition
of $R'_n$, we know that $\nu_1$ separates
$f(D(y_n,r))$ for any $r\in (R'_n,1]$.
Hence, by Lemma~\ref{lem:berkinf}.b,
for any such $r$, there is a disk
$\Dbar(z,r'')\subseteq D(y_n,r)$ such that
$f_*(\nu(z,r'')) = \nu_1$.  On the other hand,
by Lemma~\ref{lem:berkfin}, there are only finitely
many such disks $\Dbar(z,r'')$ in $D(y_n,(R'_n + 1)/2)\subsetneq D(0,1)$.
If no such disks were contained in $\Dbar(y_n,R'_n)$, then letting
$r>R'_n$ be the minimum distance from such a disk to $y_n$,
the conclusion of Lemma~\ref{lem:berkinf}.b would fail
for $D(y_n,r)$, contradicting the fact that
$\nu_1$ separates $f(D(y_n,r))$.
Thus, our claim is proved, and the desired disk
$\Dbar(z_{n+1},R''_{n+1})$ exists.

Because $f_*(\nu(z_{n+1},R''_{n+1}))=\nu_1$,
we have $G(\nu(z_{n+1},R''_{n+1}))\geq \mu^2$.
By our assumption that $G(\nu(y_n,R'_n))< \mu^2$,
we have
$\Dbar(z_{n+1},R''_{n+1})\subsetneq \Dbar(y_n,R'_n)$;
in particular, $R''_{n+1}<R'_n$.

Define $g(r) = G(\nu(z_{n+1},r))$, so that $g$
is continuous on the interval $[R''_{n+1},R'_n]$, with
$g(R''_{n+1}) \geq \mu^2 > g(R'_n)$.  
Therefore, $g$ has an absolute maximum at some radius
$R_{n+1}\in [R''_{n+1},R'_n)$.
Without loss, we may assume
$R_{n+1}$ is the largest radius in $[R''_{n+1},R'_n]$ for which
$g$ attains its maximum.  Thus, $g$ is strictly decreasing
on $[R_{n+1},R_{n+1} + \eps]$ for some $\eps>0$; by Lemma~\ref{lem:monom},
we have
$$N_0(z_{n+1},R_{n+1}) + N_{\alpha}(z_{n+1},R_{n+1})
+ N_1(z_{n+1},R_{n+1}) + N_{\infty}(z_{n+1},R_{n+1})
  > 2 N_{\text{ram}}(z_{n+1},R_{n+1}).$$
By Lemma~\ref{lem:center}, there is a point
$y_{n+1}\in D(x,R)$ such that for every
$r\in [0,R_{n+1}]$,
$$N_0(y_{n+1},r) + N_{\alpha}(y_{n+1},r)
+ N_1(y_{n+1},r) + N_{\infty}(y_{n+1},r)
  > 2 N_{\text{ram}}(y_{n+1},r).$$
By Lemma~\ref{lem:berkext},
$\nu_1$ separates $f(D(z_{n+1},R''_{n+1}+\eps))$
for every $\eps>0$.
Since 
$\Dbar(z_{n+1},R''_{n+1})\subseteq\Dbar(y_{n+1},R_{n+1})$,
it follows that $\nu_1$ separates $f(D(y_{n+1},R_{n+1}+\eps))$
for every $\eps>0$.
Thus, the pair $(y_{n+1},R_{n+1})$ satisfies
properties~(i) and~(iii) above, so that our inductive
process may repeat.

To show that the process must eventually end,
we first claim that $y_n\not\in\Dbar(y_{n+1},R'_{n+1})$
for every $n\geq 1$.
Otherwise, because $y_{n+1}\in\Dbar(y_n,R'_n)$
and $R'_{n+1}\leq R_{n+1}<R'_n$, we would have
$D(y_n,R'_n) = D(y_{n+1},R'_n)$.  However,
$\nu_1$ separates $f(D(y_{n+1},R'_n))$,
by condition~(i) for $(y_{n+1},R'_{n+1})$.
At the same time, $\nu_1$ does not
separate $f(D(y_n,R'_n))$, by condition~(ii)
for $(y_n,R'_n)$.  This contradiction proves the claim.
In particular, all of the $\{y_n\}$ are distinct.

Observe that by choosing $r$ sufficiently small
in property~(iii), each $y_n$ must have
$f(y_n)\in\{0,\alpha,1,\infty\}$.
By Lemma~\ref{lem:berkfin}, there are only
finitely many such points $y_n$ in $\Dbar(0,R)$;
hence, the process must eventually stop.
Thus, we obtain a pair $(x,R')$ satisfying properties
(i)--(iv) of Lemma~\ref{lem:mudisk}.  By the same Lemma,
then, we are done.
\qed

\begin{cor}
\label{cor:ahlfors}
Let $K$, $k$, $p$, and $E_p$ be as in Theorem~\ref{thm:ahlfors}.
Let $U_1,U_2,U_3,U_4 \subseteq \PK$ be
four pairwise disjoint open disks.
Let $\nu_1\in\PKBerk$ such that no connected
component of $\PKBerk \setminus \{ \nu_1\}$ intersects
more than two of $U_1,U_2,U_3,U_4$.
Define $C_2$ as in Theorem~\ref{thm:ahlfors}.

Let $f$ be a nonconstant
meromorphic function on $K$ such that
for any point $\nu\in\PKBerk\setminus\{\infty\}$
for which $f_*(\nu)=\nu_1$, we have $L(f,\nu) \geq C_2$.
Then there is an open disk $U\subseteq D(0,1)$ such that $f$ is one-to-one
on $U$ and $f(U)=U_i$ for some $i=1,2,3,4$.
\end{cor}

{\bf Proof.}
If $\charact K = 0$,
then the hypothesis that $f$ is nonconstant implies that there is
some $x\in K$ such that $f^{\#}(x)>0$.

On the other hand, if $\charact K = p >0$, then there are many
nonconstant functions for which $f'(z)=0$; any function which
can be written as $f(z)=g(z^p)$ with $g$ meromorphic has this property.
To avoid this situation, we first claim that $\nu_1$ separates $f(D(0,r))$
for some $r>0$.

If not, change coordinates so that $f(0)=0$ and $\nu_1=\nu(0,\rho)$
for some $\rho>0$.  Then $f$ is holomorphic on $K$ with image
contained in $D(0,\rho)$.  Because $f$ is nonconstant,
Lemma~\ref{lem:monom} implies that $\|f\|_{\nu(0,r)}\geq cr^n$
for some $c>0$ and $n\geq 1$.  Thus, we may choose $r>0$ large enough
that $\|f\|_{\nu(0,r)}> \rho$, which is a contradiction
and proves our claim.

By Lemma~\ref{lem:berkinf}.b, there is some $\nu\in\DBerk(0,r)$
such that $f_*(\nu)=\nu_1$.  Since $\charact K>0$, we have
$C_2>0$; by hypothesis, then, $L(f,\nu)> 0$, so that $f^{\#}(\nu)> 0$.
Therefore there is a point $x\in K$ such that $f^{\#}(x)>0$.

In any characteristic, then, given the point $x$ above,
define $C_1$ as in Theorem~\ref{thm:ahlfors}.
By an affine change of coordinates (moving $x$ to $0$
and scaling appropriately), we may assume that
$f^{\#}(0)>C_1$.  The result then follows by restricting $f$ to
the open unit disk and invoking Theorem~\ref{thm:ahlfors}.
\qed

\section{Examples}
\label{sect:ex}

Theorem~\ref{thm:ahlfors} differs from its complex counterpart
in several noticeable ways.  First and foremost, only four islands
are required, as opposed to the five in the complex case.
As observed in \cite{Ben7}, Example~6, it would be impossible
to reduce the number further, to three islands.  On the other
hand, in the case of positive residue characteristic, the
non-archimedean theorem requires the extra condition
that $L(f,\nu)\geq C_2$ for any $\nu$ mapping to $\nu_1$.
A similar condition, that $L(f,\nu(0,r))\geq C_2$ for {\em some} $r$,
is required for the holomorphic version \cite{Ben7}.  In this
section, we present examples to illustrate both that the lower bound
of $C_2$ is essentially sharp, and that
it is not enough to assume only one $\nu$
mapping to $\nu_1$ satisfies the inequality, even
if the number of islands is increased.

For the sharpness of the constant $C_2$,
in Example~\ref{ex:ellip} we will consider only the case
that $\charact k = p > 0=\charact K$; the map in question
is very much analogous to that of \cite{Ben7}, Example~5.  As in \cite{Ben7},
we show that the bounds given in Theorem~\ref{thm:ahlfors} are sharp
if $p=2$ and almost sharp (except possibly for the constant
$|p|^{-E_p}$) if $p\geq 3$.
If $\charact k=0$, then the lower bound $C_2$ is vacuously sharp.
We conjecture that examples analogous to \cite{Ben7},
Examples~3 and~4,
would prove the sharpness of $C_2$ in the cases that
$\charact K=p>0$.

\begin{example}
\label{ex:ellip}
Suppose $\charact K=0$ by $\charact k=p\geq 2$.  (For example,
suppose $K=\Cp$, the completion of an algebraic closure
of the $p$-adic rationals $\Qp$.)  Let $E$ be an elliptic
curve defined over $K$ with identity point $O$, and let $n\geq 1$
be an integer.  Assume that $E$ has good ordinary reduction.
(The conclusions we will reach also hold for multiplicative
reduction, but that case is slightly more complicated because
there are a number of different points $\nu$ which will map to
the point $\nu_1$ we will choose shortly.)
Let $E_1$ be the of points which map to $\bar{O}$ under
reduction.

If we identify $E_1$ with the open unit disk, then
by the characteristics of $K$ and $K$, and
because the corresponding formal group has height~1, there are
$p$-torsion points $\{P_1,\ldots,P_{p-1}\}$
at distance $|p|^{1/(p-1)}$ from $O$ in $E_1$.  Moreover, there
are no nontrivial torsion points closer than $|p|^{(p-1)}$ to $O$.

The multiplication-by-$n$ map $[n]:E\rightarrow E$
has the property that $[n](-P) = -[n](P)$ for any point $P$ on $E$.
Meanwhile, the group $\{\pm 1\}$ acts on the curve $E$ (with
$-1$ taking $P$ to $-P$) with quotient $\PP^1$.  It follows
that there is a map
$f_n:\PP^1\rightarrow\PP^1$ for which
$$
\begin{CD}
E  @>[n]>> E \\
@VhVV	 @VhVV \\
\PP^1 @>f_n>> \PP^1
\end{CD}
$$
commutes, where $h$ is quotient map.
The function $f_n$ is known as a Latt\`{e}s map to dynamicists.
It is a rational function of (geometric) degree $n^2$.


Let $a_1,a_2,a_3,a_4$ be the images under $h$ of the $2$-torsion
points $E[2]$ of $E$.  Note that $E[2]$ is the set of
ramification points of $h$.
For convenience, choose coordinates on $\PP^1$ so that $h(O)=a_1=0$
and so that $h(E_1)=D(0,1)$.  Let $\nu_1=\nu(0,1)$; then because
$E$ has good reduction, we have $(f_n)_*(\nu)=\nu_1$ if and only
if $\nu=\nu_1$.  Furthermore, an examination of the map
$[n]$ restricted to the formal group $E_1$ shows that
$L(f_n,\nu_1)=|n|$.
Let $\mu=|n|$.

Let $C_n=E[2n]\setminus E[2]$ be the set of $2n$-torsion points
which are not $2$-torsion, and let $B_n\subseteq \PK$ be the
image of $C_n$ under $h$.  By considering the ramification of
the map $h$, and knowing that $f_n$ must have exactly
$2\deg f_n -2 = 2n^2-2$ critical points, it is easy to check
that $B_n$ is the set of critical points of $f_n$,
and that each point of $B_n$ maps $2$-to-$1$ to its image.
One can also check that $f_n'(0)=n^2$.

Because $h$ maps $E_1$ two-to-one onto $D(0,1)$,
it is not difficult to show that $|h(P_i)|=|p|^{2/(p-1)}$.
Similarly, the lack of nontrivial torsion points
in $E_1$ at distance less than $|p|^{1/(p-1)}$ from $O$
implies that
\begin{equation}
\label{eq:notors}
\left( \{a_2,a_3,a_4\}\cup B_n \right) \cap D(0,|p|^{2/(p-1)}) = \emptyset .
\end{equation}
In fact, if $p\nmid n$ and $p\geq 3$, then
$\left( \{a_2,a_3,a_4\}\cup B_n \right) \cap D(0,1) = \emptyset$.


Meanwhile, $B_n\cup\{a_1,a_2,a_3,a_4\}$ is precisely the preimage
of $\{a_1,a_2,a_3,a_4\}$ under $f_n$.
Thus, if we choose $U_1,U_2,U_3,U_4$ to be disjoint
disks containing $a_1,a_2,a_3,a_4$ in Theorem~\ref{thm:ahlfors},
then the only disks $U$ that could map one-to-one onto any $U_i$
would have to contain exactly one of $a_1,a_2,a_3,a_4$ and cannot
intersect $B_n$.  By translating on $E$ by $2$-torsion points,
it suffices to consider only one-to-one mappings of $f_n$ from a
disk $U$ containing $a_1=0$ to $U_1\subseteq D(0,1)$.  In fact, the
preimage disk $U$ we should consider is the largest disk about $0$
which contains no points in $B_n\cup\{a_2,a_3,a_4\}$.
If $p|n$ or $p=2$, then
by equation~\eqref{eq:notors}, this disk is $U=D(0,|p|^{2/(p-1)})$.
By Lemma~\ref{lem:holom}, the largest possible one-to-one image
disk is $U_1=D(0,|n|^2\cdot |p|^{2/(p-1)})$, since $|f'_n(0)|=|n|^2$.
On the other hand, if $p\nmid n$ and $p\geq 3$,
the domain disk is $U=D(0,1)$, and its image is $U_1=D(0,1)$.

If $p=2$, then there is one nontrivial $2$-torsion point $P_1$ in $E_1$.
Its image $h(P_1)=a_2$ is the point $\alpha$ of Section~\ref{sect:thm}.
In addition, no component of $\PK\setminus\{\nu_1\}$ contains more
than two of $a_1,a_2,a_3,a_4$; thus, $f_n$ satisfies the conditions
of Corollary~\ref{cor:ahlfors}.  Since $|\alpha|=|p|^{2/(p-1)}$,
the radius of the image disk $U_1$ described in the previous
paragraph is exactly $|\alpha|\cdot \mu^2$.  Thus, the Ahlfors
radius is $s=\mu^2$, which is exactly the
lower bound in Lemma~\ref{lem:mudisk}.

If $p\geq 3$, we have $|\alpha|=1$, since none
of $\{a_2,a_3,a_4\}$ lie in $D(0,1)$.
If $p\nmid n$, then $\mu=1$, so that the lower bound
from Lemma~\ref{lem:mudisk} for the
radius $s$ of the image disk is $1$, which
is exactly the radius of the disk $U_1$ found above.
Finally, if $p|n$, then the image disk has radius
$s=\mu^2 |p|^{2/(p-1)}$, which is only slightly larger
than the Lemma~\ref{lem:mudisk} lower bound of $\mu^2 |p|^{2E_p}$.
\end{example}

Our final example will illustrate that having only
one $\nu$ for which $f_*(\nu)=\nu_1$ with $L(f,\nu)$ bounded below
by some fixed amount is
not enough to guarantee an islands theorem, regardless of how
many islands there are, how small they are, or how small the lower
bound on $L(f,\nu)$ is.
As in Example~\ref{ex:ellip}, Example~\ref{ex:warp} is only
for the case that $\charact k >\charact K=0$.  Of course,
if $\charact k=0$, then the condition that $L(f,\nu)\geq 0$ is vacuous,
as previously noted, so there will be no counterexamples, as we
already know the theorem is already true without hypothesis~(b)
in that case.
On the other hand, if $\charact K=p>0$,
then we imagine that examples similar to the
following one may be constructed.

\begin{example}
\label{ex:warp}
Assume that $\charact k = p > 0 = \charact K$.
For any integer $N\geq 0$ and any radius $0<s\leq 1$, we select $N+2$ islands
as follows.  Set $a_0=0\in K$.
For each $i=1,\ldots, N$, choose $a_i\in K$ with $|a_i|=1$
and, for each $i\neq j$, $|a_i-a_j|=1$.
The first $N+1$ islands will be the open disks $D(a_i,s)$, with
$i=0,\ldots, N$.  The final island will be
$\PK\setminus\Dbar(0,1/s)$, which is the open disk of spherical
radius $s$ centered at $\infty$.  Let $\nu_1=\nu(0,1)$,
which separates each island from every other.

Pick $n\geq 1$ large enough so that $|p^n|<s$, and choose
$b\in K$ so that $0<|b|<s$.  Let $c=-b^{1+p^n}$.

For each $i=1,\ldots, N$, define
$$f_i(z) = \frac{(z-a_i)^{p^n} + c}{(z-a_i)^{p^n}}
  = 1 + \frac{c}{(z-a_i)^{p^n}}.$$
Then, define
$$f_0(z) = \frac{z^{1+p^n} + c}{z^{p^n}} = z + \frac{c}{z^{p^n}},
\quad\text{and}\quad
f(z) = \prod_{i=0}^N f_i(z).$$
The reader may check that $f_*(\nu_1)=\nu_1$ and
$L(f,\nu(0,1))=1$; recall from
Lemma~\ref{lem:Gbound} that this is the maximum value
$L$ could ever attain.
%

The only preimages of $\infty$ are $a_0,\ldots,a_N$,
all of which are critical points.  Thus,
the island at $\infty$ has no one-to-one preimages.
In addition, if $|z-a_i|\geq 1$ for all $i=0,\ldots, N$, then
it is easy to see that $|f(z)-a_i|\geq 1$ also.  In particular,
any preimages of the islands $D(a_i,s)$ must lie in the disks
$D(a_j,1)$.

Next, observe that there are $1+p^n$ preimages of $0$ in $D(0,1)$,
namely the roots of $f_0(z)= 0$, which are all of the form
$\zeta^j b$, where $\zeta$ is a primitive $(1+p^n)$-root of unity.
In particular, the largest open disk about any such root which
contains no other such roots has radius $|b|$.
It is easy to compute that $|f'(\zeta^j b)| = 1$, and therefore
the image of that largest open disk is a disk of radius $|b|<s$.
Thus, the island at $0$ is not a one-to-one image of a disk inside
$D(0,1)$.

Similarly, for any fixed $i=1,\ldots, N$, there are $p^n$ preimages
of $0$ in $D(a_i,1)$, namely the roots of $(z-a_i)^{p^n}=-c$.
Those roots are of the form $x=a_i + \omega^j d$, where $\omega$
is a $p^n$-root of unity, and $d$ is a $p^n$-root of $-c$.
Any $\omega^j$ is distance $|p|^{1/(p-1)}$ from the nearest other
$\omega^{\ell}$, so that the largest disk containing exactly
one root of $(z-a_i)^{p^n}=-c$ has radius $|p|^{1/(p-1)}|c|^{1/p^n}$.
We can compute that $|f'(x)| = |p|^n\cdot |c|^{-1/p^n}$, so that
the largest possible one-to-one image disk has radius
$|p|^{n + (1/(p-1))} < |p|^n < s$.  Thus, there are no one-to-one
preimages of the island at $0$ anywhere in $K$.

For $1\leq i\leq N$, the preimages of $a_i$ in $D(0,1)$ must be points
$x\in D(0,1)$ satisfying $|f_0(x)|=1$.  Since $|x|<1$, this must
mean $|x|=|c|^{1/p^n}$.  Because of the pole at $0$, the largest
open disk about $x$ which could conceivably map onto the island
$D(a_i,s)$ must have radius at most $|c|^{1/p^n}$.  (In fact, it
will have slightly smaller radius than that, but the bound
of $|c|^{1/p^n}$ will suffice for our purposes.)
Some computation using the above value for $|x|$ shows that
$|f'(x)| \leq \max \{1, |p|^n |c|^{-1/p^n} \}$, 
and therefore the largest possible one-to-one image disk
has radius
$$\max \left\{|c|^{1/p^n}, |p|^n\right\} \leq
\max \left\{|b|, |p|^n\right\} < s,$$
which fails to cover the island.

Before considering the final case, observe that
$$f(z) - z = z\cdot \left[ 1 - \left( 1 + \frac{c}{z^{1+p^n}} \right)
\prod_{j=1}^N \left( 1 + \frac{c}{(z-a_j)^{p^n}} \right) \right],
$$
so that
$$
|f(z) - z| \leq |z| \cdot \max \left\{
\left| \frac{c}{z^{1+p^n}} \right|, 
\left| \frac{c}{(z-a_1)^{p^n}} \right|, \ldots
\left| \frac{c}{(z-a_N)^{p^n}} \right|\right\}.
$$
It follows that if $|c|^{1/(1+p^n)}<|x-a_j| < 1$ for some $j=1,\ldots, N$
and some $x\in K$, then $|f(x)-x| < |x-a_j|$.  In particular,
no such $x$ can have image $f(x)$ in any of the $N+2$ islands.

We are now ready to consider the final possibility, that there
is a preimage of an island $D(a_i,s)$ in the disk $D(a_j,1)$
for some $i,j\in\{1,\ldots, N\}$, where $i$ and $j$ may or may not
be equal.
In such a case, we have a point $x$ in $D(x,a_j)$ with $f(x)=a_i$.
This means that $|f_j(x)|=1$, so that $|x-a_j|\geq |c|^{1/p^n}$.
From this bound, it follows that
$|f'(x)| \leq \max \{1, |p^n|/|x-a_j|\}$.
Meanwhile, the largest possible disk mapping onto $D(a_i,s)$
cannot contain the pole at $a_j$, so that its radius must be
at most $|x-a_j|$.  Thus, the largest possible one-to-one
image disk about $a_i$ has radius at most
$$|f'(x)| \cdot |x-a_j| \leq \max \{|x-a_j|,|p^n|\}
\leq \max \{|c|^{1/(1+p^n)},|p^n|\}
= \max \{|b|, |p^n|\} <s,$$
where the second inequality is by the previous paragraph.
Thus, none of the $N$ islands has a one-to-one preimage anywhere
in $K$, in spite of the fact that $f_*(\nu_1)=\nu_1$
with $L(f,\nu(0,1))=1$.

\end{example}

An examination of the proof of Theorem~\ref{thm:ahlfors} applied
to Example~\ref{ex:warp} reveals why the hypothesis that {\em every}
(or at least many) 
$\nu$ mapping to $\nu_1$ must have $L(f,\nu)\geq \mu$.
For that choice of $f$,
we start from $R=1$ (or whatever smaller radius $\Dbar(0,1)$ is moved
to after $f$ is scaled as described in the proof of
Corollary~\ref{cor:ahlfors})
and move inward, searching for a disk satisfying conditions (i)--(iv)
of Lemma~\ref{lem:mudisk}.  Properties (i), (iii),
and (iv) already apply to $\Dbar(0,R)$, but property (ii) does not.
Lemma~\ref{lem:center} would select the new center $y_1$ to be one
of the roots of $f=a_i$ in $D(a_i,1)$ for some $i=0,\ldots, N$;
the inductive process would begin with $R_1=1$ and $y_1$ being one
such root.  The minimal radius $R'_1$ would be the smallest radius
about $y_1$ for which $f(\Dbar(y_1,R'_1))$ contained points outside
$D(a_i,1)$.  That is, $R'_1=|y_1-a_i|$, which we saw to be at
most $|b|$.  Even though the inequality of property (iii) holds,
the degree of $r$ in $G$ is positive (in fact, equal to $1$) for
$r\in [R'_1,R_1]$, so that as $r$ shrinks from $R_1$ down to $R'_1$,
$G$ also shrinks from $1$ down to $R'_1\leq |b|$.  Thus, although
$f$ is at last one-to-one on $D(y_1,R'_1)$, the image is too small
because there are poles too close to $y_1$, and the value of $G$
(as well as $L$, along with it) has shrunk too much.  There are
type~II points $\nu$ which separate the smaller disk $D(y_1,R'_1)$,
but without hypothesis (b.) of the Theorem, their $G$-values are
too small.  Thus, if we try to shrink to a disk $D(z_2,R''_2)$
according to the algorithm, we have no guarantee that $G$ increases,
and therefore we have no guarantee that property (iii) holds.

On the other hand,
a modified version of the same
example, with $f_0(z) = z + c/z^{p^n}$ replaced by $f_0(z)=z+c/z^{p^n-1}$,
has the same pathology of points $\nu$ mapping to $\nu_1$ with
small $G(\nu)$ in each of the disks $D(a_i,1)$ for
$i=1,\ldots, N$.  This time, however, for $i=0$, there is
one extra $\nu$ mapping to $\nu_1$
in $D(0,1)$ which {\em does} satisfy $G(\nu)=1$.  As a result, the conclusion
of Theorem~\ref{thm:ahlfors} holds for the modified example,
because there are disks in $D(0,1)$ which map one-to-one onto,
say, $D(a_1,1)$.  Thus, it is conceivable that some condition weaker
than 
``$L(f,\nu)\geq C_2$ for {\em all} $\nu$ mapping to $\nu_1$''
but stronger than
``$L(f,\nu)\geq C_2$ for {\em some} $\nu$ mapping to $\nu_1$''
would suffice.  We leave the existence of such a condition
as an open question.

\bibliographystyle{plain}

\end{document}